\documentclass[11pt,a4paper]{article}

\usepackage{mathtools}
\usepackage{authblk} % author
\usepackage{algpseudocode,algorithmicx,algorithm}

\usepackage{mathrsfs}
\usepackage{latexsym,bm}
\usepackage{amsmath,amsfonts,amsmath,amssymb,amsthm}
\usepackage{extarrows}

\usepackage{graphicx,subfigure,epstopdf,float}
\usepackage{enumerate,cases,multirow}
\usepackage{makecell}
\usepackage{caption}

\usepackage{longtable,colortbl,arydshln,threeparttable}
\definecolor{mygray}{gray}{.9}

\usepackage{indentfirst}
\usepackage[top=25mm,bottom=20mm,left=25mm,right=20mm]{geometry}
\usepackage{cite}

\usepackage{listings}

\usepackage{makeidx}        % generate an index, automatically
\usepackage{booktabs}
\usepackage[bookmarksnumbered,colorlinks,citecolor=red,linkcolor=red,hyperindex,linktocpage=true]{hyperref}

\newcommand{\bb}{\boldsymbol}

\def \e {\mathrm{e}}

\newcounter{parentalgorithm}

\makeatother

\newtheorem{example}{\bf Example}[section]
\theoremstyle{remark}
\newtheorem{remark}{\bf Remark}[section]

\numberwithin{equation}{section}

\usepackage{algorithm}
\usepackage{listings,xcolor} % ³ÌÐò´úÂë
\usepackage{textcomp}
\usepackage[framed,numbered,autolinebreaks,useliterate]{mcode} % mcode.sty is needed

\begin{document}
%\begin{CJK*}{GBK}{song}

\title{\bf varFEM: variational formulation based programming for finite element methods in Matlab}
\author{Yue Yu\thanks{terenceyuyue@sjtu.edu.cn}}
\affil{School of Mathematical Sciences, Institute of Natural Sciences, MOE-LSC, Shanghai Jiao Tong University, Shanghai, 200240, P. R. China.}
\date{}
\maketitle

\begin{abstract}
  This paper summarizes the development of varFEM, which provides a realization of the programming style in FreeFEM by using the Matlab language.
\end{abstract}

%\textbf{Keywords.} Virtual element method, Polygonal meshes, Three dimensions, Matlab

%\renewcommand{\contentsname}{\centering \sanhao \hei \quad }
%\tableofcontents

\section{Introduction}

FreeFEM is a popular 2D and 3D partial differential equations (PDE) solver based on finite element methods (FEMs) \cite{FreeFEM}, which has been used by thousands of researchers across the world. The highlight is that the programming language is consistent with the variational formulation of the underlying PDEs, referred to as the variational formulation based programming in this article. We intend to develop an FEM package in a similar way of FreeFEM using the language of Matlab, named varFEM. The similarity here only refers to the programming style of the main or test script, not to the internal architecture of the software.

This programming paradigm is usually organized in an object-oriented language, which makes it difficult for readers or users to understand and modify the code, and further redevelop the package (although it is a good way to develop softwares). Upon rethinking the process of finite element programming, it becomes clear that the assembly of the stiffness matrix and load vector essentially reduces to the numerical integration of some
typical bilinear and linear forms with respect to basis functions. In this regard, the package $i$FEM, written in Matlab, has provided robust, efficient, and easy-following codes for many mathematical and physical problems in both two and three dimensions \cite{ChenL-iFEM-2009}.
On this basis, we successfully developed the variational formulation based programming for the conforming $\mathbb{P}_k$-Lagrange ($k\le 3$) FEMs in two dimensions by utilizing the Matlab language. The underlying idea can be generalized to other types of finite elements for both two- and three-dimensional problems on unstructured simplicial meshes. The package is accessible on \url{https://github.com/Terenceyuyue/varFEM} (see the varFEM folder).

The article is organized as follows. In Section \ref{sect:varFEM}, we introduce the basic idea in varFEM through a model problem.  In Section \ref{sect:tutorial}, we demonstrate the use of varFEM for several typical examples, including a complete implementation of the model problem, the vector finite element for the linear elasticity, the mixed FEMs for bihamonic equation and Stokes problem, and the iterative scheme for the heat equation. We also demonstrate the ability of varFEM to solve complex problems in Section \ref{sect:freefem}.

\section{Variational formulation based programming in varFEM} \label{sect:varFEM}

We introduce the variational formulation based programming in varFEM via a model problem, so as to facilitate the underlying design idea.

\subsection{Programming for a model problem}
Let $\Omega = (0,1)^2$ and consider the second-order elliptic problem:
\begin{equation}\label{modelP}
\begin{cases}
   - \nabla  \cdot (a\nabla u) + cu = f \quad & \text{in}~~\Omega ,  \\
  u = g_D\quad & \text{on}~~\Gamma _D,  \\
  g_Ru + a\partial _nu = g_N\quad  & \text{on}~~\Gamma _R,
\end{cases}
\end{equation}
where $\Gamma_D$ and $\Gamma _R = \partial \Omega \backslash \Gamma_D$ are the Dirichlet boundary and Robin boundary, respectively.
For brevity, we refer to $g_R$ as the Robin boundary data function and $g_N$ as the Neumann boundary data function.
For homogenous Dirichlet boundary condition, the variational problem is to find $u \in V: = H_0^1(\Omega )$ such that
\[a(v,u) = \ell(v),\quad v \in V,\]
where
\[a(v,u) = \int_\Omega  {a\nabla v \cdot \nabla u} {\rm d}\sigma  + \int_\Omega  {cvu} {\rm d}\sigma  + \int_{\Gamma _R}g_Rvu {\rm d}s,\]
\[\ell (v) = \int_\Omega  fv {\rm d}\sigma  + \int_{\Gamma _R} g_N v {\rm d}s.\]
Here, the test function is placed in the first entry of $a(\cdot,\cdot)$ since
\[a(v,u) = \bb{v}^T\bb{Au},\quad \bb{A} = (a(\varphi _i,\phi _j))_{m \times n},\]
where $\bb{v}= (v_1,\cdots,v_m)^T$ and $\bb{u} = (u_1,\cdots,u_n)^T$, with
\[v = \sum\limits_{i = 1}^m v_i\varphi _i,\quad u = \sum\limits_{i = 1}^n u_i\phi _i.\]

\subsubsection{The assembly of bilinear forms}

The first step is to obtain the stiffness matrix associated with the bilinear form
\[\int_\Omega  a\nabla v \cdot \nabla u {\rm d}\sigma  + \int_\Omega  cvu {\rm d}\sigma \]
on the approximated domain, where for simplicity we have used the original notation $\Omega$ to represent a triangular mesh.
The computation in varFEM reads
\vspace{-0.8cm}
\begin{lstlisting}
% Omega
Coef  = {a, c};
Test  = {'v.grad', 'v.val'};
Trial = {'u.grad', 'u.val'};
kk = int2d(Th,Coef,Test,Trial,Vh,quadOrder);
\end{lstlisting}
Here, \mcode{Th} represents the triangular mesh, which provides some necessary auxiliary data structures. We set up the triple \mcode{(Coef,Test,Trial)}, for the coefficients, test functions and trial functions in variational form, respectively. It is obvious that \mcode{v.grad} is for $\nabla v$ and \mcode{v.val} is for $v$ itself. The routine \mcode{int2d.m}  computes the stiffness matrix corresponding to the bilinear form on the two-dimensional region, i.e.
\[A = (a_{ij}),\quad a_{ij} = a(\Phi _i,\Phi _j),\]
where $\Phi _i$ are the global shape functions of the finite element space \mcode{Vh}. The integral of the bilinear form, as
$(\nabla \Phi _i, \nabla \Phi _j)_\Omega$, will be approximated by using the Gaussian quadrature formula with \mcode{quadOrder} being the order of accuracy.

The second step is to compute the stiffness matrix for the bilinear form on the Robin boundary $\Gamma_R$:
\[\int_{\Gamma_R} g_Rvu {\rm d}s.\]
The code can be written as follows.
\vspace{-0.8cm}
\begin{lstlisting}
% Gamma_R
Th.elem1d = Th.bdEdgeType{1};
Th.elem1dIdx = Th.bdEdgeIdxType{1};
Coef  = {g_R};
Test  = {'v.val'};
Trial = {'u.val'};
kk = kk + int1d(Th,Coef,Test,Trial,Vh,quadOrder);
\end{lstlisting}
Here, \mcode{int1d.m} gives the contribution to the stiffness matrix on the one-dimensional boundary edges of the mesh.
Note that we must provide the connectivity list \mcode{elem1d} of the boundary edges of $\Gamma_R$ and the associated indices \mcode{elem1dIdx} in the data structure \mcode{edge} introduced later.

\subsubsection{The assembly of linear forms}

For the linear forms, we first consider the integral for the source term:
\[\int_\Omega  fv {\rm d}\sigma.\]
 The load vector can be assembled as
\vspace{-0.8cm}
\begin{lstlisting}
% Omega
Coef = pde.f;  Test = 'v.val';
ff = int2d(Th,Coef,Test,[],Vh,quadOrder);
\end{lstlisting}
We set \mcode{Trial = []} to indicate the linear form.

The computation of the load vector associated with the Neumann boundary data function $g_N$, i.e.,
\[\int_{\Gamma_R} g_Nv {\rm d}s\]
reads
\vspace{-0.8cm}
\begin{lstlisting}
% Gamma_R
Coef = g_N;  Test = 'v.val';
ff = ff + int1d(Th,Coef,Test,[],Vh,quadOrder);
\end{lstlisting}

\subsection{Data structures for triangular meshes}

We adopt the data structures given in $i$FEM \cite{ChenL-iFEM-2009}. All related data are stored in the Matlab structure \mcode{Th}, which is computed by using the subroutine \mbox{FeMesh2d.m} as
\vspace{-0.8cm}
\begin{lstlisting}
Th = FeMesh2d(node,elem,bdStr);
\end{lstlisting}

The triangular meshes are represented by two basic data structures \mcode{node} and \mcode{elem}, where \mcode{node} is an $\mcode{N} \times 2$ matrix with the first and second columns contain $x$- and $y$-coordinates of the nodes in the mesh, and \mcode{elem} is an $\mcode{NT} \times 3$ matrix recording the vertex indices of each element in a counterclockwise order, where \mcode{N} and \mcode{NT} are the numbers of the vertices and triangular elements.

In the current version, we only consider the $\mathbb{P}_k$-Lagrange finite element spaces with $k$ up to 3. In this case, there are two important data structures \mcode{edge} and \mcode{elem2edge}. In the matrix \mcode{edge(1:NE,1:2)}, the first and second rows contain indices of the starting and ending points. The column is sorted in the way that for the $k$-th edge, \mcode{edge(k,1)<edge(k,2)} for $k=1,2,\cdots,\mcode{NE}$. The matrix \mcode{elem2edge} establishes the map of local index of edges in each triangle to its global index in matrix \mcode{edge}. By convention, we label three edges of a triangle such that the $i$-th edge is opposite to the $i$-th vertex. We refer the reader to \url{https://www.math.uci.edu/~chenlong/ifemdoc/mesh/auxstructuredoc.html} for some detailed information.

To deal with boundary integrals, we first exact the boundary edges from \mcode{edge} and store them in matrix \mcode{bdEdge}.
In the input of \mcode{FeMesh2d}, the string \mbox{bdStr} is used to indicate the interested boundary part in \mcode{bdEdge}. For example, for the unit square $\Omega = (0,1)^2$,
\begin{itemize}
  \item \mcode{bdStr = 'x==1'} divides \mcode{bdEdge} into two parts: \mcode{bdEdgeType\{1\}}
gives the boundary edges on $x=1$, and \mcode{bdEdgeType\{2\}} stores the remaining part.
  \item \mcode{bdStr = \{'x==1','y==0'\} } separates the boundary data \mcode{bdEdge} into three parts: \mcode{bdEdgeType\{1\}} and \mcode{bdEdgeType\{2\}} give the boundary edges on $x=1$ and $y=0$, respectively, and \mcode{bdEdgeType\{3\}} stores the remaining part.
  \item \mcode{bdStr = []} implies that \mcode{bdEdgeType\{1\} = bdEdge}.
\end{itemize}

We also use \mcode{bdEdgeIdxType} to record the index in matrix \mcode{edge}, and \mcode{bdNodeIdxType} to store the nodes for respective boundary parts. Note that we determine the boundary of interest by the coordinates of the midpoint of the edge, so \mcode{'x==1'} can also be replaced by a statement like \mcode{'x>0.99'}.

\subsection{Code design of \mcode{int2d.m} and \mcode{assem2d.m}}

In this article we only discuss the implementation of the bilinear forms in two dimensions.

\subsubsection{The scalar case: \mcode{assem2d.m}} \label{subsect:assem2d}

In this subsection we introduce the details of writing the subroutine \mcode{assem2d.m} to assemble a two-dimensional scalar bilinear form
\[a(v,u),\quad v = \varphi _i, ~~u = \phi _j,\quad i = 1, \cdots ,m,~~j = 1, \cdots ,n,\]
where the test function $v$ and the trial function $u$ are allowed to match different finite element spaces, which can be found in mixed finite element methods for Stokes problems.
For the scalar case, \mcode{assem2d.m} is essentially the same as \mcode{int2d.m}, while the later one can be used to deal with vector cases like linear elasticity problems. To handle different spaces, we write \mcode{Vh = \{'P1', 'P2' \}} for the input of \mcode{assem2d.m}, where \mcode{Vh\{1\}} is for $v$ and \mcode{Vh\{2\}} is for $u$. For simplicity, it is also allowed to write \mcode{Vh = 'P1'} when $v$ and $u$ are in the same space.

Let us discuss the case where $v$ and $u$ lie in the same space. Suppose that the bilinear form contains only first-order derivatives. Then the possible combinations are
\[\int_K {avu} {\rm d}x,\quad \int_K av_xu {\rm d}x,\quad \int_K av_yu {\rm d}\sigma,\]
\[\int_K avu_x {\rm d}\sigma,\quad \int_K av_xu_x {\rm d}\sigma,\quad \int_K av_yu_x {\rm d}\sigma,\]
\[\int_K avu_y {\rm d}\sigma,\quad \int_K av_xu_y {\rm d}\sigma,\quad \int_K av_yu_y {\rm d}\sigma.\]
Of course, we often encounter the gradient form
\[\int_K  a\nabla v \cdot \nabla u {\rm d}\sigma = \int_K  a(v_xu_x + v_yu_y) {\rm d}\sigma.\]

We take the second bilinear form as an example. Let
\[a_K(v,u) = \int_K a v_xu {\rm d}\sigma,\]
and consider the $\mathbb{P}_1$-Lagrange finite element. Denote the local basis functions to be $\phi _1,\phi _2,\phi _3$. Then the local stiffness matrix is
\[A_K = \int_K a  \begin{bmatrix}
  \partial _x\phi _1 \\
  \partial _x\phi _2 \\
  \partial _x\phi _3
\end{bmatrix} \begin{bmatrix}\phi _1 & \phi _2 & \phi _3 \end{bmatrix} {\rm d}\sigma.\]
Let
\[v_1 = \partial _x\phi _1, ~~v_2 = \partial _x\phi _2,~~v_3 = \partial _x\phi _3;\quad
u_1 = \phi _1,~~u_2 = \phi _2,~~u_3 = \phi _3.\]
Then
\[A_K = (k_{ij})_{3 \times 3}, \quad k_{ij} = \int_K av_iu_j {\rm d}\sigma.\]
The integral will be approximated by the Gaussian quadrature rule:
\[
k_{ij} = \int_Kav_iu_j {\rm d}\sigma = | K |\sum\limits_{p = 1}^{n_g} w_pa(x_p,y_p)v_i(x_p,y_p)u_j(x_p,y_p) ,
\]
where $(x_p,y_p)$ is the $p$-th quadrature point. In the implementation, we in advance store the quadrature weights and the values of basis functions or their derivatives in the following form:
\[w_p,\qquad \mcode{vi}(:,p) = \begin{bmatrix}
  v_i|_{(x_p^1,y_p^1)} \\
  v_i|_{(x_p^2,y_p^2)} \\
   \vdots  \\
  v_i|_{(x_p^{\text{NT}},y_p^{\text{NT}})}
\end{bmatrix},\qquad \mcode{uj}(:,p) =  \begin{bmatrix}
  u_j|_{(x_p^1,y_p^1)} \\
  u_j|_{(x_p^2,y_p^2)} \\
   \vdots  \\
  u_j|_{(x_p^{\text{NT}},y_p^{\text{NT}})}
\end{bmatrix},\qquad p = 1, \cdots ,n_g,\]
where $\mcode{vi}$ associated with $v_i$ is of size $\mcode{NT} \times \mcode{ng}$ with the $p$-th column given by $\mcode{vi}(:,p)$.
Let $\mcode{weight} = [w_1, w_2, \cdots, w_{n_g}]$ and \mcode{ww = repmat(weight, NT, 1)}. Then $k_{ij}$ for $a=1$ can be computed as
\vspace{-0.8cm}
\begin{lstlisting}
k11 = sum(ww.*v1.*u1,2);
k12 = sum(ww.*v1.*u2,2);
k13 = sum(ww.*v1.*u3,2);
k21 = sum(ww.*v2.*u1,2);
k22 = sum(ww.*v2.*u2,2);
k23 = sum(ww.*v2.*u3,2);
k31 = sum(ww.*v3.*u1,2);
k32 = sum(ww.*v3.*u2,2);
k33 = sum(ww.*v3.*u3,2);
K = [k11,k12,k13, k21,k22,k23, k31,k32,k33];
\end{lstlisting}
Here we have stored the local stiffness matrix $A_K$ in the form of $[k_{11},k_{12},k_{13},k_{21},k_{22},k_{23},k_{31},k_{32},k_{33}]$, and stacked the results of all cells together. By adding the contribution of the area, one has
\vspace{-0.8cm}
\begin{lstlisting}
Ndof = 3;
K = repmat(area,1,Ndof^2).*K;
\end{lstlisting}

For the variable coefficient case, such as $a(x,y) = x+y$,  one can further introduce the \textbf{coefficient matrix} as
\vspace{-0.8cm}
\begin{lstlisting}
cf = @(pz) pz(:,1) + pz(:,2); % x+y;
cc = zeros(NT,ng);
for p = 1:ng
    pz = lambda(p,1)*z1 + lambda(p,2)*z2 + lambda(p,3)*z3;
    cc(:,p) = cf(pz);
end
\end{lstlisting}
where \mcode{pz} are the quadrature points on all elements. The above procedure can be implemented as follows.
\vspace{-0.8cm}
\begin{lstlisting}
K = zeros(NT,Ndof^2);
s = 1;
v = {v1,v2,v3}; u = {u1,u2,u3};
for i = 1:Ndof
    for j = 1:Ndof
        vi = v{i}; uj = u{j};
        K(:,s) =  area.*sum(ww.*cc.*vi.*uj,2);
        s = s+1;
    end
end
\end{lstlisting}
The bilinear form is assembled by using the build-in function \mcode{sparse.m} as in $i$FEM. In this case, the code is given as
\vspace{-0.8cm}
\begin{lstlisting}
ss = K(:);
kk = sparse(ii,jj,ss,NNdof,NNdof);
\end{lstlisting}
Here, the triple \mcode{(ii, jj, ss)} is called the sparse index. Please refer to the following link: \url{https://www.math.uci.edu/~chenlong/ifemdoc/fem/femdoc.html}.

\begin{remark}
  For the case where $v$ and $u$ are in different spaces, one just needs to modify the basis functions and the number of local degrees of freedom accordingly. The code can be presented as
\vspace{-0.8cm}
\begin{lstlisting}
s = 1;
for i = 1:Ndofv
    for j = 1:Ndofu
        vi = vbase{i}; uj = ubase{j};
        K(:,s) =  K(:,s) + area.*sum(ww.*cc.*vi.*uj,2);
        s = s+1;
    end
end
\end{lstlisting}
\end{remark}

In varFEM, we use \mcode{Base2d.m} to load the information of \mcode{vi} and \mcode{uj}, for example, the following code gives the values of $\partial_x \phi$, where $\phi$ is a local basis function.
\vspace{-0.8cm}
\begin{lstlisting}
v = 'v.dx';
vbase = Base2d(v,node,elem,Vh{1},quadOrder); % v1.dx, v2.dx, v3.dx
\end{lstlisting}

\subsubsection{The vector case: \mcode{int2d.m}} \label{subsect:int2d}

Let us consider a typical bilinear form for linear elasticity problems, given as
\[a_K(\bb{v}, \bb{u}) := \int_K  \bb{\varepsilon}(\boldsymbol v) : \bb{\varepsilon }(\boldsymbol u) \text{d}\sigma,\]
where $\bb{v} = (v_1,v_2)^T$, $\bb{u} = (u_1,u_2)^T$, and
\begin{align}
\bb{\varepsilon}(\boldsymbol v) : \bb{\varepsilon }(\boldsymbol u)
& = v_{1,x}u_{1,x} + v_{2,y}u_{2,y} + \frac{1}{2}(v_{1,y} + v_{2,x})(u_{1,y} + u_{2,x}) \label{originalForm}\\
& = v_{1,x}u_{1,x} + v_{2,y}u_{2,y} + \frac{1}{2}(v_{1,y}u_{1,y}  + v_{1,y}u_{2,x}+ v_{2,x}u_{1,y} + v_{2,x}u_{2,x}). \label{extendedForm}
\end{align}
The stiffness matrix can be assembled as
 \vspace{-0.8cm}
\begin{lstlisting}
Coef  = {1, 1, 0.5, 0.5, 0.5, 0.5};
Test  = {'v1.dx', 'v2.dy', 'v1.dy', 'v1.dy', 'v2.dx', 'v2.dx'};
Trial = {'u1.dx', 'u2.dy', 'u1.dy', 'u2.dx', 'u1.dy', 'u2.dx'};
kk = int2d(Th, Coef, Test, Trial, Vh, quadOrder);
\end{lstlisting}
We also provide the subroutine \mcode{getExtendedvarForm.m} to get the extended combinations \eqref{extendedForm} from \eqref{originalForm}, which has been included in \mcode{int2d.m}. Therefore, the bilinear form can be directly assembled as
 \vspace{-0.8cm}
\begin{lstlisting}
Coef = { 1, 1, 0.5 };
Test  = {'v1.dx', 'v2.dy', 'v1.dy + v2.dx'};
Trial = {'u1.dx', 'u2.dy', 'u1.dy + u2.dx'};
kk = int2d(Th, Coef, Test, Trial, Vh, quadOrder);
\end{lstlisting}

In the rest of this subsection, we briefly discuss the sparse index. Let $\bb v = (v_1,v_2,v_3)$ and $\bb u = (u_1,u_2,u_3)$, and suppose that $a(\bb v,\bb u)$ is a bilinear form. Note that, in general, $v_i (i=1,2,3)$ can be in different spaces, but $v_i$ and $u_i$ are in the same space, otherwise the resulting stiffness matrix is not a square matrix. The stiffness matrix after blocking has the following correspondence:
   \[a(\bb v,\bb u) \qquad \leftrightarrow \qquad [\bb v_1,\bb v_2,\bb v_3]\left[ \begin{array}{*{20}{c}}
  A_{11}&A_{12}&A_{13} \\
  A_{21}&A_{22}&A_{23} \\
  A_{31}&A_{32}&A_{33}
\end{array} \right]\left[ \begin{array}{*{20}{c}}
  \bb u_1 \\
  \bb u_2 \\
  \bb u_3
\end{array} \right],\]
  where $\bb u_i$ is the vector of degrees of freedom of $u_i$. It is easy to see that $A_{ij}$ can obtained as in scalar case by assembling all pairs that contain $(v_i,u_j)$ in $a(\bb v,\bb u)$.

  Let the sparse index for $A_{ij}$ be $(\bb{i}_{ij},\bb{j}_{ij}, \bb{s}_{ij})$. Let the numbers of rows and columns of $A_{ij}$ be $m_i$ and $n_j$, respectively. Then the final sparse assembly index \mcode{ii} and \mcode{jj} can be written in block matrix as
  \[\left[ \begin{array}{*{20}{c}}
  \bb{i}_{11}&\bb{i}_{12}&\bb{i}_{13} \\
  \bb{i}_{21} + m_1&\bb{i}_{22} + m_1&\bb{i}_{23}+ m_1 \\
  \bb{i}_{31} + m_2&\bb{i}_{32} + m_2&\bb{i}_{33} + m_2
\end{array} \right]\quad \mbox{and} \quad \left[ \begin{array}{*{20}{c}}
  \bb{j}_{11}&\bb{j}_{12} + n_1&\bb{j}_{13} + n_2 \\
  \bb{j}_{21}&\bb{j}_{22} + n_1&\bb{j}_{23} + n_2 \\
  \bb{j}_{31}&\bb{j}_{32} + n_1&\bb{j}_{33} + n_2
\end{array} \right],\]
and obtained by straightening them as a column vector along the row vectors.

\section{Tutorial examples} \label{sect:tutorial}

In this section, we present several examples to demonstrate the use of varFEM.

\subsection{Poisson-type problems}

We now provide the complete implementation of the model problem \eqref{modelP}. The function file reads
\vspace{-0.8cm}
\begin{lstlisting}
function uh = varPoisson(Th,pde,Vh,quadOrder)

%% Assemble stiffness matrix
% Omega
Coef  = {pde.a, pde.c};
Test  = {'v.grad', 'v.val'};
Trial = {'u.grad', 'u.val'};
kk = assem2d(Th,Coef,Test,Trial,Vh,quadOrder); % or assem2d

% Robin data
bdStr = Th.bdStr;
if ~isempty(bdStr)
    Th.elem1d = Th.bdEdgeType{1};
    Th.elem1dIdx = Th.bdEdgeIdxType{1};

    Coef  = pde.g_R;
    Test  = 'v.val';
    Trial = 'u.val';
    kk = kk + assem1d(Th,Coef,Test,Trial,Vh,quadOrder); % or assem1d
end

%% Assemble the right hand side
% Omega
Coef = pde.f;
Test = 'v.val';
ff = int2d(Th,Coef,Test,[],Vh,quadOrder);
% Neumann data
if ~isempty(bdStr)
    %Coef = @(p) pde.g_R(p).*pde.uexact(p) + pde.a(p).*(pde.Du(p)*n');

    fun = @(p) pde.g_R(p).*pde.uexact(p);
    Cmat1 = interpEdgeMat(fun,Th,quadOrder);
    fun = @(p) repmat(pde.a(p),1,2).*pde.Du(p);
    Cmat2 = interpEdgeMat(fun,Th,quadOrder);
    Coef = Cmat1 + Cmat2;

    ff = ff + assem1d(Th,Coef,Test,[],Vh,quadOrder);
end

%% Apply Dirichlet boundary conditions
g_D = pde.g_D;
on = 2 - 1*isempty(bdStr); % 1 for bdStr= [], 2 for bdStr = 'x==0'
uh = apply2d(on,Th,kk,ff,Vh,g_D);
\end{lstlisting}

In the above code, the structure \mcode{pde} stores the information of the PDE, including the exact solution \mcode{pde.uexact}, the gradient \mcode{pde.Du}, etc. The Neumann data function is $g_N = g_Ru + a\partial _nu$, which varies on the boundary edges. For testing purposes, we compute this function by using the exact solution. In Lines 31-35, we use the subroutine \mcode{interpEdgeMat.m} to derive the coefficient matrix as in Subsect.~\ref{subsect:assem2d}. We remark that the \mcode{Coef} has three forms:
\begin{enumerate}
  \item A function handle or a constant.
  \item The numerical degrees of freedom of a finite element function.
  \item A coefficient matrix resulting from the numerical integration.
\end{enumerate}
In the computation, the first two forms in fact will be transformed to the third one.

The boundary edges will be divided into at most two parts. For example, when \mcode{bdStr = 'x==0'}, the Robin boundary part is given by \mcode{elem1d = Th.bdEdgeType\{1\}}. The remaining part is for the Dirichlet boundary condition, with the index \mcode{on = 1} given for \mcode{Th.bdEdgeType}.

The test script is presented as follows.
 \vspace{-0.8cm}
\begin{lstlisting}
%% Parameters
maxIt = 5;
N = zeros(maxIt,1);
h = zeros(maxIt,1);
ErrL2 = zeros(maxIt,1);
ErrH1 = zeros(maxIt,1);

%% Generate an intitial mesh
[node,elem] = squaremesh([0 1 0 1],0.5);
bdStr = 'x==0';

%% Get PDE data
pde = Poissondatavar;
g_R = @(p) 1 + p(:,1) + p(:,2); % 1 + x + y
pde.g_R = g_R;

%% Finite Element
i = 1; % 1,2,3
Vh = ['P', num2str(i)];
quadOrder = i+2;
for k = 1:maxIt
    % refine mesh
    [node,elem] = uniformrefine(node,elem);
    % get the mesh information
    Th = FeMesh2d(node,elem,bdStr);
    % solve the equation
    uh = varPoisson(Th,pde,Vh,quadOrder);
    % record and plot
    N(k) = size(elem,1);
    h(k) = 1/(sqrt(size(node,1))-1);
    if N(k) < 2e3  % show mesh and solution for small size
        figure(1);
        showresult(node,elem,pde.uexact,uh);
        drawnow;
    end
    % compute error
    ErrL2(k) = varGetL2Error(Th,pde.uexact,uh,Vh,quadOrder);
    ErrH1(k) = varGetH1Error(Th,pde.Du,uh,Vh,quadOrder);
end

%% Plot convergence rates and display error table
figure(2);
showrateh(h,ErrH1,ErrL2);
fprintf('\n');
disp('Table: Error')
colname = {'#Dof','h','||u-u_h||','|u-u_h|_1'};
disptable(colname,N,[],h,'%0.3e',ErrL2,'%0.5e',ErrH1,'%0.5e');
\end{lstlisting}

In the \mcode{for} loop, we first load or generate the mesh, which immediately returns the matrix \mcode{node} and \mcode{elem} to the Matlab workspace. Then we set up the boundary conditions to get the structural information. The subroutine \mcode{varPoisson.m} is the function file containing all source code to implement the FEM as given before. When obtaining the numerical solutions, we can visualize the solutions by using the subroutines \mcode{showresult.m}. We then calculate
the discrete $L^2$ and $H^1$ errors via the subroutines \mcode{varGetL2Error.m} and \mcode{varGetH1Error.m}. The procedure is completed by verifying the rate of convergence through \mcode{showrateh.m}.

\begin{figure}[H]
  \centering
  \includegraphics[scale=0.5]{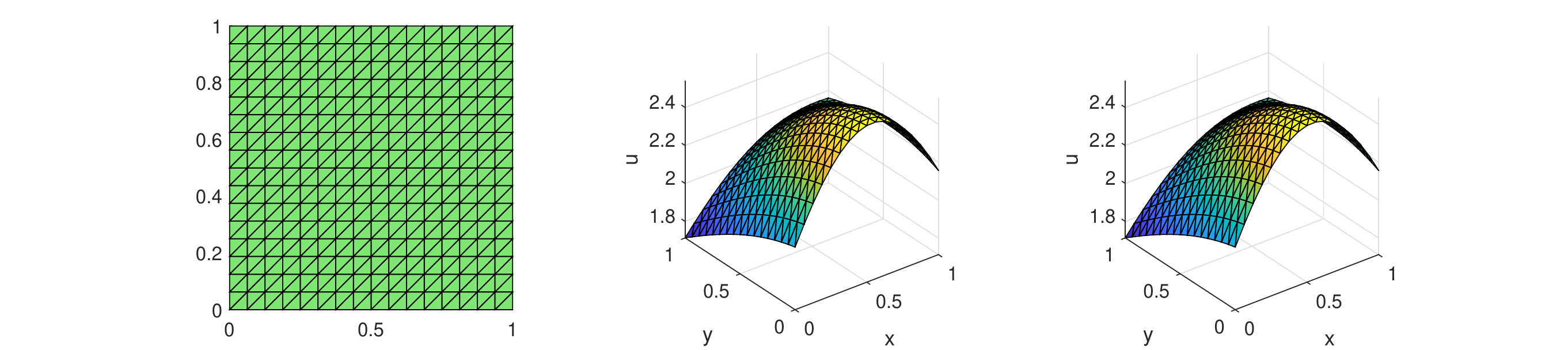}\\
  \caption{The nodal values of exact and numerical solutions for the model problem}\label{fig:Poisson_values}
\end{figure}

\begin{figure}[!htb]
  \centering
  \subfigure[$i=1$]{\includegraphics[scale=0.35]{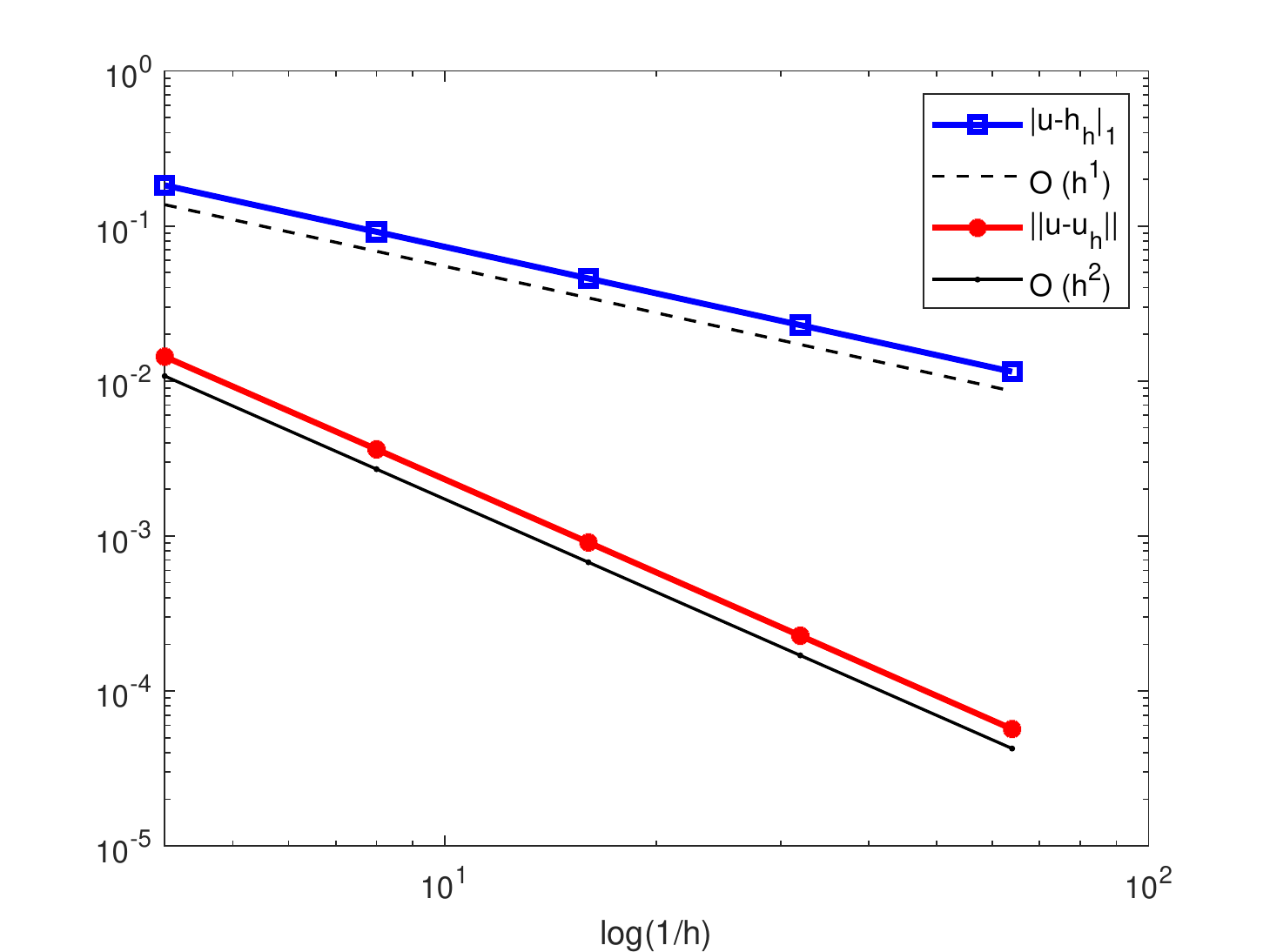}}
  \subfigure[$i=2$]{\includegraphics[scale=0.35]{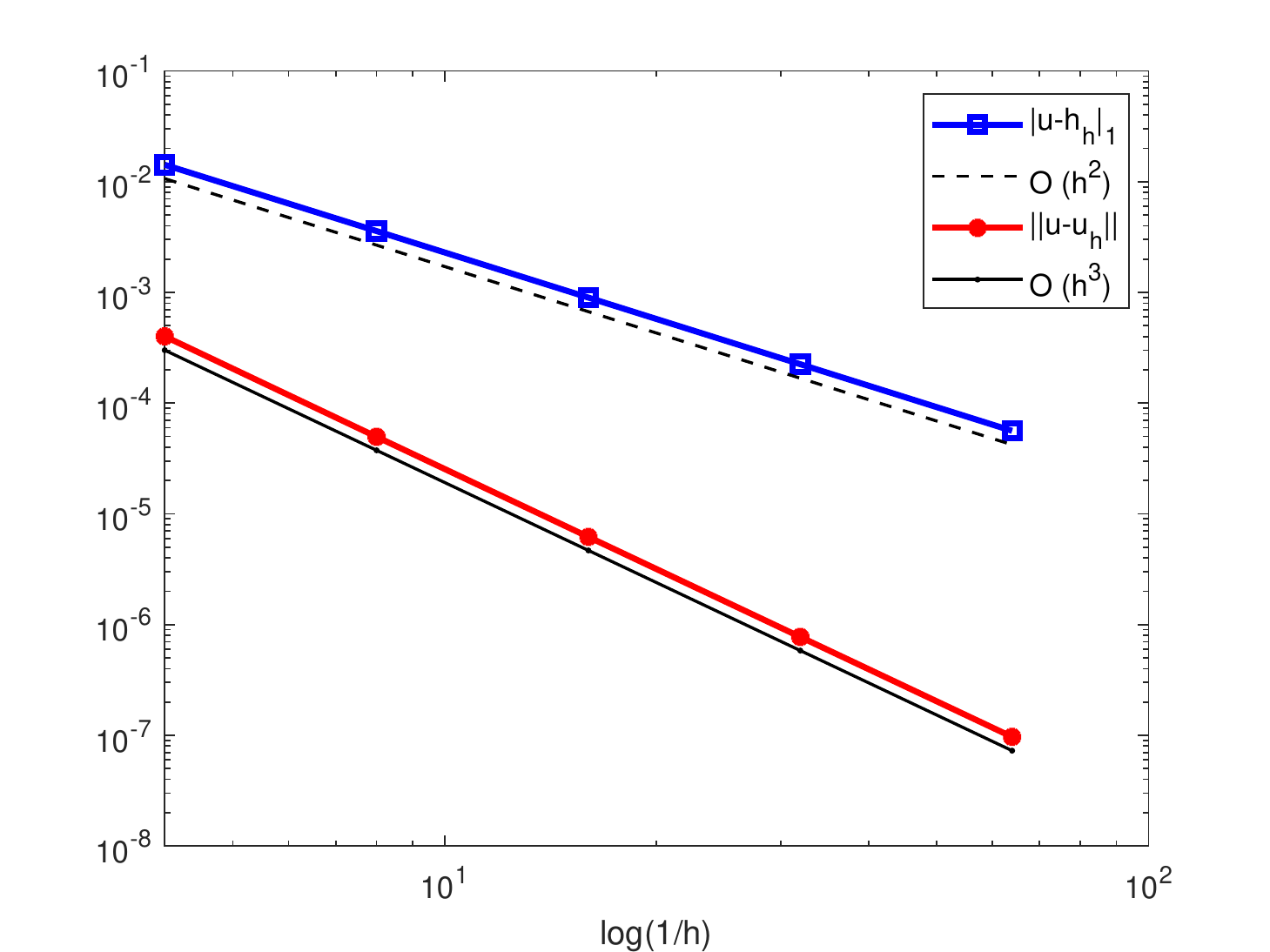}}
  \subfigure[$i=3$]{\includegraphics[scale=0.35]{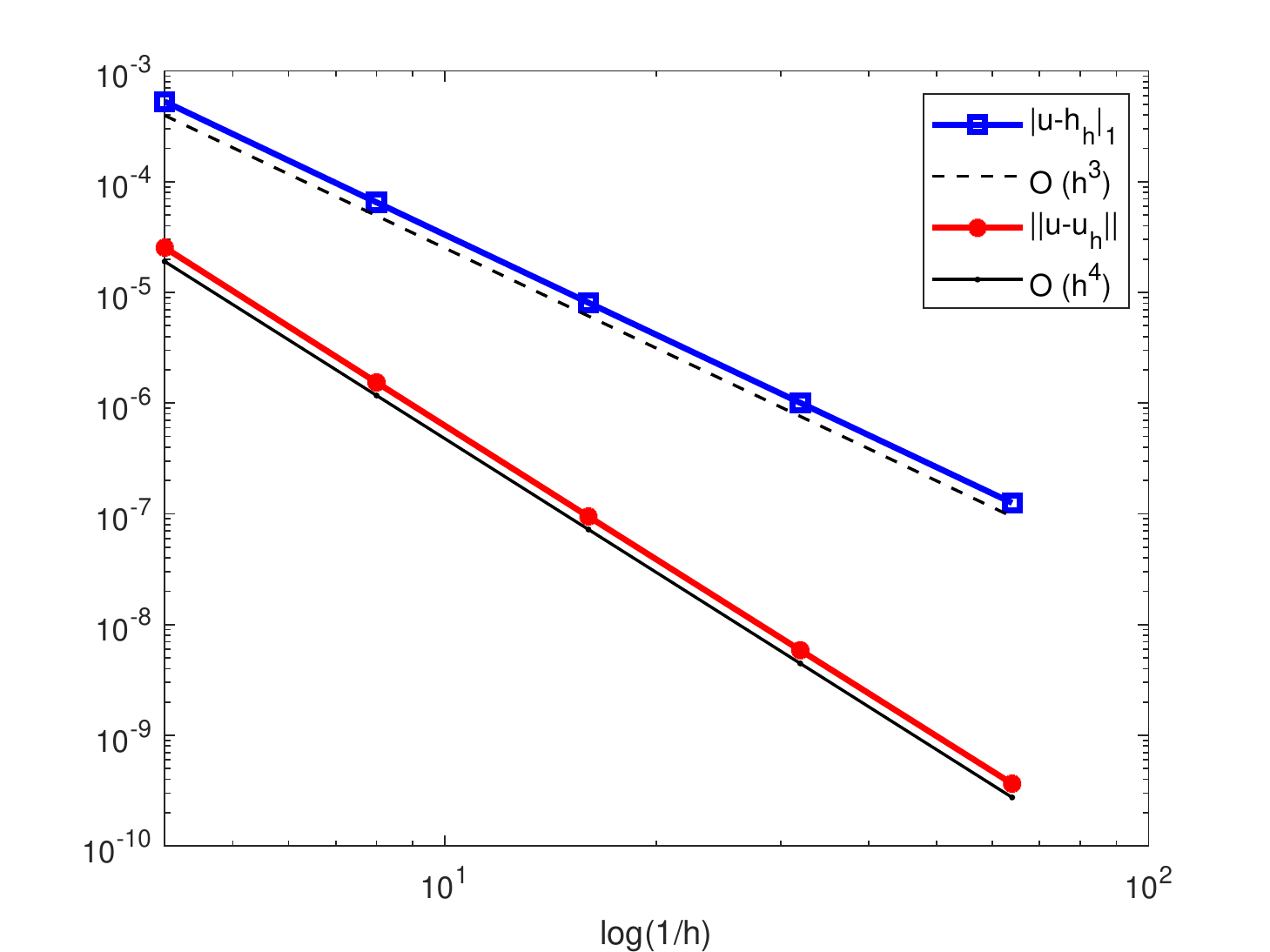}}\\
  \caption{Convergence rates for the model problem.}\label{fig:Poisson_rate}
\end{figure}

The test script can be easily used to compute $\mathbb{P}_i$-Lagrange element for $i=1,2,3$ (see Line 18 in the test script). The nodal values for the model problem are displayed in Fig.~\ref{fig:Poisson_values}. The rates of convergence are shown in Fig.~\ref{fig:Poisson_rate}, from which we observe the optimal convergence for all cases.

We also print the errors on the Matlab command window, with the results for $\mathbb{P}_3$ element given below:
{\scriptsize\begin{verbatim}
        Table: Error
            #Dof        h         ||u-u_h||      |u-u_h|_1
            ____    _________    ___________    ___________

              32    2.500e-01    2.53816e-05    5.28954e-04
             128    1.250e-01    1.53744e-06    6.50683e-05
             512    6.250e-02    9.46298e-08    8.07922e-06
            2048    3.125e-02    5.87001e-09    1.00695e-06
            8192    1.562e-02    3.65444e-10    1.25698e-07
\end{verbatim}}

We next consider the example with a circular domain or an L-shaped domain. Such a domain can be generated by using the pdetool as
\vspace{-0.8cm}
\begin{lstlisting}
%% Generate an intitial mesh
g = 'circleg'; %'lshapeg';
[p,e,t] = initmesh(g,'hmax',0.5);
node = p'; elem = t(1:3,:)';
bdStr = 'x>0'; % string for Neumann
\end{lstlisting}

\begin{figure}[H]
  \centering
  {\includegraphics[scale=0.5]{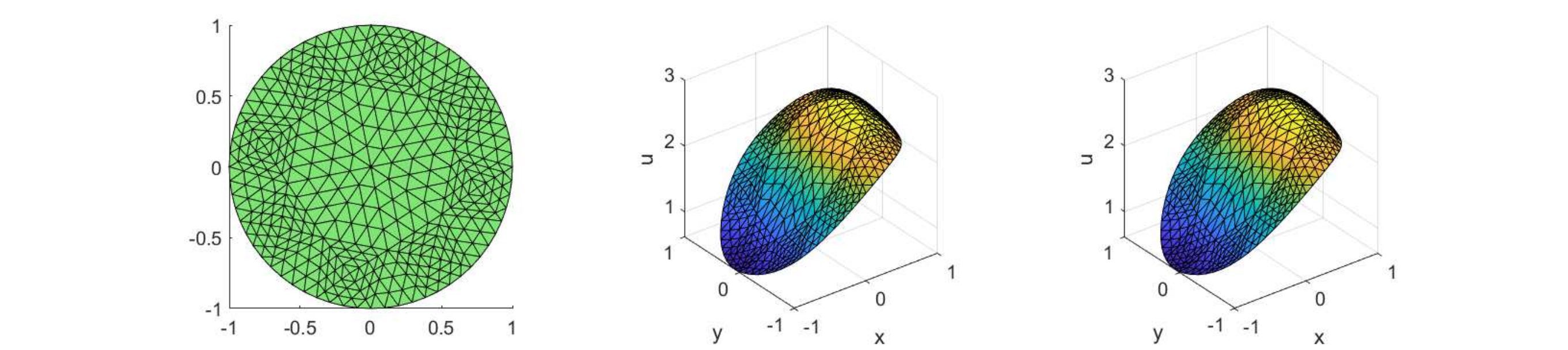}}\\
  {\includegraphics[scale=0.5]{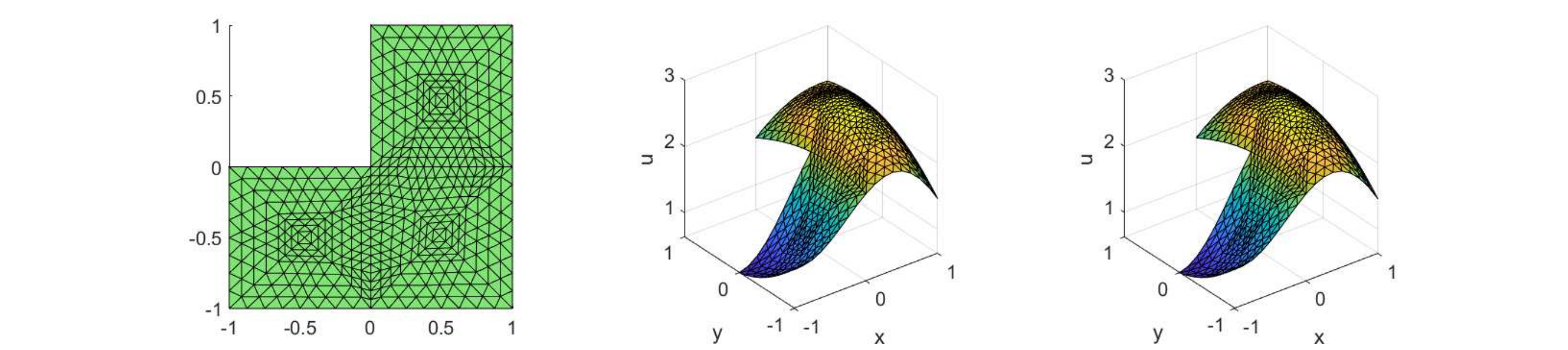}}\\
  \caption{The nodal values of exact and numerical solutions on circular domain or L-shaped region.}\label{fig:Poisson_values_circle}
\end{figure}

The results are given in Fig.~\ref{fig:Poisson_values_circle}

\subsection{Linear elasticity problems}

The linear elasticity problem is
\begin{equation}\label{linearElasticity}
\begin{cases}
   - {\rm div} \bb{\sigma } = \bb{f} \quad & \text{in}~~~\Omega ,  \\
  \bb{u} = {\bb{0}} \quad & \text{on}~~~\Gamma_0,  \\
  \bb{\sigma n} = \bb{g} \quad & \text{on}~~~\Gamma_1,
\end{cases}
\end{equation}
where $\bb{n} = (n_1,n_2)^T$ denotes the outer unit vector normal to $\partial\Omega$. The constitutive relation for linear elasticity is
\[\bb{\sigma}(\bb{u}) = 2\mu\bb{\varepsilon}(\bb{u}) + \lambda({\rm div} \bb{u})\bb{I},\]
where $\bb{\sigma} = (\sigma_{ij})$ and $\bb{\varepsilon} = (\varepsilon_{ij})$ are the second order stress and strain tensors, respectively, satisfying $\varepsilon_{ij} = \frac{1}{2}(\partial_i u_j + \partial_j u_i)$, $\lambda$ and $\mu$ are the Lam\'{e} constants, $\bb{I}$ is the identity matrix, and ${\rm div}\bb{u} = \partial_1 u_1 + \partial_2 u_2$.

\subsubsection{The programming in scalar form}

The vector problem can be solved in block form by using \mcode{assem2d.m} as scalar cases. The equilibrium equation in \eqref{linearElasticity} can also be written in the form
\begin{equation}\label{displacementType}
 - \mu \Delta \bb{u} - (\lambda  + \mu ){\text{grad}}({\rm div} \bb{u}) = \bb{f} \quad {\text{in}}~~\Omega,
\end{equation}
which is referred to as the displacement type in what follows. In this case, we only consider $\Gamma_0 = \Gamma := \partial \Omega$. The first term $\Delta \bb{u}$ can be treated as the vector case of the Poisson equation. The variational formulation is
\[
  \mu \int_\Omega  \nabla \bb u \cdot \nabla \bb v {\rm d}x + (\lambda  + \mu )\int_\Omega  (\text{div}~\bb u)(\text{div}~\bb v) {\rm d}\sigma
   = \int_\Omega  \bb{f} \cdot \bb v {\rm d}\sigma.
\]

The first term of the bilinear form can be split into
\[\mu \int_\Omega \nabla u_1 \cdot \nabla v_1 \text{d}\sigma\quad \mbox{and} \quad
\mu \int_\Omega \nabla u_2 \cdot \nabla v_2 \text{d}\sigma. \]
They generate the same matrix, denoted $A$, corresponding to the blocks  $A_{11}$ and $A_{22}$, respectively. The computation reads
\vspace{-0.8cm}
\begin{lstlisting}
% (v1.grad, u1.grad), (v2.grad, u2.grad)
cf = 1;
Coef  = cf;  Test  = 'v.grad';  Trial = 'u.grad';
A = assem2d(Th,Coef,Test,Trial,Vh,quadOrder);
\end{lstlisting}

The second term of the bilinear form has the following combinations:
\[\int_\Omega  v_{1,x}u_{1,x} \text{d}x,\quad
\int_\Omega  v_{1,x}u_{2,y} \text{d}x,\quad
\int_\Omega  v_{2,y}u_{1,x} \text{d}x\quad
 \int_\Omega  v_{2,y}u_{2,y} \text{d}x,\]
which correspond to $A_{11}$, $A_{12}$, $A_{21}$ and $A_{22}$, respectively, and can be computed as follows.
\vspace{-0.8cm}
\begin{lstlisting}
% (v1.dx, u1.dx)
cf = 1;
Coef  = cf;  Test  = 'v.dx';  Trial = 'u.dx';
B1 = assem2d(Th,Coef,Test,Trial,Vh,quadOrder);
% (v1.dx, u2.dy)
cf = 1;
Coef  = cf;  Test  = 'v.dx';  Trial = 'u.dy';
B2 = assem2d(Th,Coef,Test,Trial,Vh,quadOrder);
% (v2.dy, u1.dx)
cf = 1;
Coef  = cf;  Test  = 'v.dy';  Trial = 'u.dx';
B3 = assem2d(Th,Coef,Test,Trial,Vh,quadOrder);
% (v2.dy, u2.dy)
cf = 1;
Coef  = cf;  Test  = 'v.dy';  Trial = 'u.dy';
B4 = assem2d(Th,Coef,Test,Trial,Vh,quadOrder);
\end{lstlisting}

The block matrix is then given by
\vspace{-0.8cm}
\begin{lstlisting}
% kk
kk = [  mu*A+(lambda+mu)*B1,         (lambda+mu)*B2;
             (lambda+mu)*B3,    mu*A+(lambda+mu)*B4   ];
kk = sparse(kk);
\end{lstlisting}

The right-hand side has two components:
\[\int_\Omega  f_1v_1 \text{d}x\quad \mbox{and}\quad \int_\Omega  f_2v_2 \text{d}x.\]
The load vector is assembled in the following way:
\vspace{-0.8cm}
\begin{lstlisting}
%% Assemble right hand side
% F1
trf = eye(2);
Coef = @(pz) pde.f(pz)*trf(:, 1);  Test = 'v.val';  % f = [f1, f2]
F1 = assem2d(Th,Coef,Test,[],Vh,quadOrder);
% F2
trf = eye(2);
Coef = @(pz) pde.f(pz)*trf(:, 2);  Test = 'v.val';
F2 = assem2d(Th,Coef,Test,[],Vh,quadOrder);
% F
ff = [F1; F2];
\end{lstlisting}

For the vector problem, we impose the Dirichlet boundary value conditions as follows.
\vspace{-0.8cm}
\begin{lstlisting}
%% Apply Dirichlet boundary conditions
g_D = pde.g_D;
on = 1;
g_D1 = @(p) g_D(p)*[1;0];
g_D2 = @(p) g_D(p)*[0;1];
gBc = {g_D1,g_D2};
Vhvec = {Vh,Vh};
u = apply2d(on,Th,kk,ff,Vhvec,gBc); % note Vhvec
\end{lstlisting}
Here, \mcode{g\_D1} is for $u_1$ and \mcode{g\_D2} is for $u_2$. Note that the finite element spaces \mcode{Vhvec} must be given in the same structure of \mcode{gBc}.

\begin{figure}[H]
  \centering
  \includegraphics[scale=0.5]{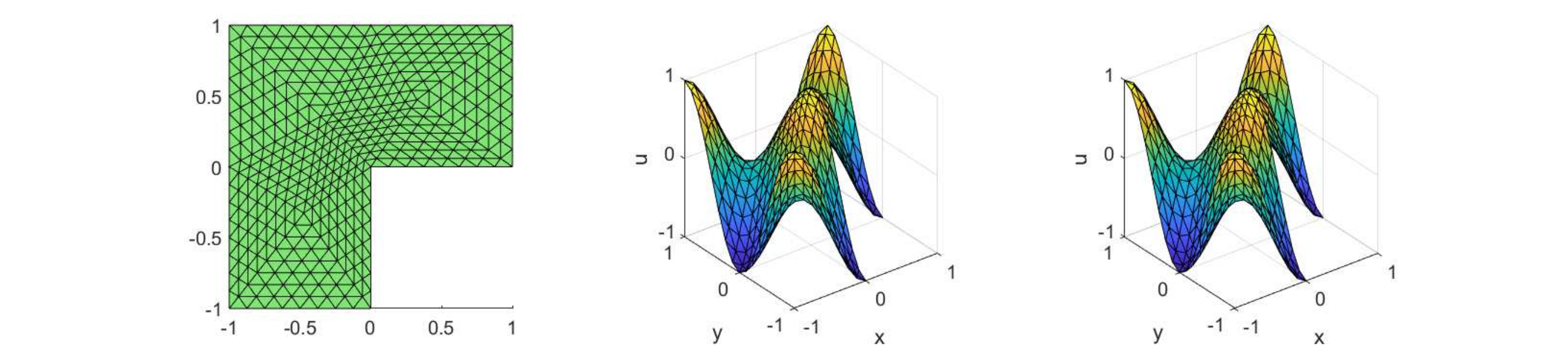}\\
  \caption{The nodal values of exact and numerical solutions for the elasticity problem of displacement type.}\label{fig:Elasticity_values1}
\end{figure}

The solutions are displayed in Fig.~\ref{fig:Elasticity_values1}.

\subsubsection{The programming in vector form}

The bilinear form is
\[
  a(\bb v,\bb u) = 2\mu \int_\Omega \varepsilon _{ij}(\bb v)\varepsilon _{ij}(\bb u) {\rm d}\sigma + \lambda \int_\Omega (\partial _iv_i)(\partial _ju_j) {\rm d}\sigma,
\]
where the summation is ommited. The computation of the first term has been given in Subsect.~\ref{subsect:int2d}, i.e.,
\vspace{-0.8cm}
\begin{lstlisting}
% (Eij(u):Eij(v))
Coef = { 1, 1, 0.5 };
Test  = {'v1.dx', 'v2.dy', 'v1.dy + v2.dx'};
Trial = {'u1.dx', 'u2.dy', 'u1.dy + u2.dx'};
A = int2d(Th,Coef,Test,Trial,Vh,quadOrder);
A = 2*mu*A;
\end{lstlisting}
The second term can be computed as
\vspace{-0.8cm}
\begin{lstlisting}
% (div u,div v)
Coef = 1;
Test  = 'v1.dx + v2.dy' ;
Trial = 'u1.dx + u2.dy';
B = int2d(Th,Coef,Test,Trial,Vh,quadOrder);
B = lambda*B;
% stiffness matrix
kk = A + B;
\end{lstlisting}

The linear form is
\[\ell (\boldsymbol v) = \int_\Omega  \bb{f} \cdot \bb v {\rm d}\sigma + \int_{\Gamma _1} \boldsymbol{g} \cdot \bb v {\rm d}s.\]
For the first term, one has
\vspace{-0.8cm}
\begin{lstlisting}
Coef = pde.f;  Test = 'v.val';
ff = int2d(Th,Coef,Test,[],Vh,quadOrder);
\end{lstlisting}
Note that we have added the implementation for $\bb{f} \cdot \bb v$ by just setting \mcode{Test = 'v.val'}.
For the second term, we first determine the one-dimensional edges:
\vspace{-0.8cm}
\begin{lstlisting}
%% Get 1D mesh for boundary integrals
bdStr = Th.bdStr;
if ~isempty(bdStr)
    Th.elem1d = Th.bdEdgeType{1};
    Th.elem1dIdx = Th.bdEdgeIdxType{1};
end
\end{lstlisting}
The coefficient matrix for the boundary integral can be computed using \mcode{interpEdgeMat.m} and the Neumann condition is then realized as
\vspace{-0.8cm}
\begin{lstlisting}
%% Assemble Neumann boundary conditions
if ~isempty(bdStr)
    g_N = pde.g_N; trg = eye(3);

    g1 = @(p) g_N(p)*trg(:,[1,3]);   Cmat1 = interpEdgeMat(g1,Th,quadOrder);
    g2 = @(p) g_N(p)*trg(:,[3,2]);   Cmat2 = interpEdgeMat(g2,Th,quadOrder);

    Coef = {Cmat1, Cmat2};  Test = 'v.val';
    ff = ff + int1d(Th,Coef,Test,[],Vh,quadOrder);
end
\end{lstlisting}
Here, the data \mcode{g_N} is stored as $g_N = [\sigma_{11}, \sigma_{22}, \sigma_{12}]$ in the structure \mcode{pde}.

The Dirichlet condition can be handled as the displacement type:
\vspace{-0.8cm}
\begin{lstlisting}
%% Apply Dirichlet boundary conditions
on = 2 - 1*isempty(bdStr);
g_D1 = @(pz) pde.g_D(pz)*[1;0];
g_D2 = @(pz) pde.g_D(pz)*[0;1];
g_D = {g_D1, g_D2};
u = apply2d(on,Th,kk,ff,Vh,g_D);
\end{lstlisting}
Note that \mcode{Vh} is of vector form.

In the test script, we set \mcode{bdStr = 'y==0 \| x==1'}. The errors for the $\mathbb{P}_3$ element are listed in Tab.~\ref{tab:ElasticityTensor}.

\begin{table}[H]
 \centering
 \caption{The $L^2$ and $H^1$ errors for the elasticity problem of tensor form ($\mathbb{P}_3$ element)} \label{tab:ElasticityTensor}
\begin{tabular}{rcccccccccc}
  \hline
  $\sharp {\rm Dof}$   &     $h$    &   ErrL2   &  ErrH1 \\
  \hline
    32   &  2.500e-01  &  6.82177e-04  &  1.34259e-02 \\
     128 &   1.250e-01 &   3.90894e-05 &   1.61503e-03 \\
     512 &   6.250e-02 &   2.32444e-06 &   1.97527e-04 \\
    2048 &   3.125e-02 &   1.42085e-07 &   2.44506e-05 \\
    8192 &   1.562e-02 &   8.79363e-09 &   3.04301e-06 \\
  \hline
\end{tabular}
\end{table}

\subsection{Mixed FEMs for the biharmonic equation}

\subsubsection{The programming in scalar form}

For the mixed finite element methods, we first consider the biharmonic equation with Dirichlet boundary conditions:
\[ \begin{cases}
  \Delta ^2u = f  \quad & {\text{in}}~~\Omega  \subset \mathbb{R}^2,  \\
  u = \partial_nu = 0 \quad & \text{on}~~\partial \Omega .
\end{cases} \]
By introducing a new variable $w = -\Delta u$, the above problem can be written in a mixed form as
\[ \begin{cases}
   - \Delta u = w,  \\
   -\Delta w = f,  \\
  u = \partial_nu = 0 \quad \text{on}~~\partial \Omega .
\end{cases} \]
The associated variational problem is: Find $(w,u) \in H^1(\Omega ) \times {H_0^1}(\Omega )=: V \times U$ such that
\begin{equation}\label{biharmix}
 \begin{cases}
 \displaystyle \int_\Omega  \nabla u \cdot \nabla \phi  {\rm d}\sigma = \int_\Omega  w\phi  {\rm d}\sigma,\quad & \phi  \in H^1(\Omega ),  \\[2mm]
 \displaystyle \int_\Omega  \nabla w \cdot \nabla \psi  {\rm d}\sigma = \int_\Omega  f\psi  {\rm d}\sigma,\quad & \psi  \in H_0^1(\Omega ).
\end{cases}
\end{equation}
Let
\[a(w,\phi ) =  - \int_\Omega  w\phi  {\rm d}\sigma \quad \mbox{and} \quad  b(\phi ,u) = \int_\Omega  \nabla \phi  \cdot \nabla u {\rm d}\sigma.\]
One has
\[ \begin{cases}
  a(w,\phi ) + b(\phi ,u)  = 0, \quad  & \phi  \in H^1(\Omega )=V,  \\
    b(w,\psi )   = (f,\psi ), \quad & \psi  \in H_0^1(\Omega )=U,
\end{cases} \]
where $a( \cdot , \cdot ): V \times V \to \mathbb{R}$ and $b( \cdot , \cdot ): V \times U \to \mathbb{R}$.

The functions $u$ and $w$ will be approximated by $\mathbb{P}_1$-Lagrange elements.  Let $N$ be the vector of global basis functions. One easily gets
\[\begin{cases}
  \bb \phi ^TA \bb{w} + \bb \phi ^TB \bb u = 0,  \\
   \bb \psi ^TB^T \bb{w} = \bb \psi ^T \bb{f},
\end{cases} \]
where
\[A = -\int_\Omega  N^TN {\rm d}\sigma, \qquad B = \int_\Omega  \nabla N^T \cdot \nabla N {\rm d}\sigma, \qquad
\bb f = \int_\Omega  N^Tf {\rm d}\sigma.\]
The system can be written in block matrix form as
\[\begin{bmatrix}
  A   & B \\
  B^T & O
\end{bmatrix}  \begin{bmatrix}
  \bb{w} \\
  \bb{u}
\end{bmatrix}  = \begin{bmatrix}
  \bb{0} \\
  \bb{f}
\end{bmatrix}.\]

\begin{remark}
 If $\partial _nu $ does not vanish on $\partial \Omega$, then the mixed variational formulation is
\[ \begin{cases}
 \displaystyle a(w,\phi ) + b(\phi ,u) = \int_{\partial \Omega} \partial_nu \phi {\rm d}s,\quad \phi  \in V,  \\
   b(w,\psi ) = (f,\psi ),\quad \psi  \in U.
\end{cases} \]
Here, $\partial _nu $ corresponds to the Neumann boundary data for $u$ in the first equation.
\end{remark}

In block form, the stiffness matrix can be computed as follows.
\vspace{-0.8cm}
\begin{lstlisting}
%% Assemble stiffness matrix
% matrix A
Coef = 1;  Test = 'v.val';  Trial = 'u.val';
A = -assem2d(Th,Coef,Test,Trial,'P1',quadOrder);
% matrix B
Coef = 1;  Test = 'v.grad';  Trial = 'u.grad';
B = assem2d(Th,Coef,Test,Trial,'P1',quadOrder);
% kk
O = zeros(size(B));
kk = [A,  B;  B', O];
kk = sparse(kk);
\end{lstlisting}
The right-hand side is given by
\vspace{-0.8cm}
\begin{lstlisting}
%% Assemble right-hand side
Coef = pde.f;  Test = 'v.val';
ff = assem2d(Th,Coef,Test,[],'P1',quadOrder);
O = zeros(size(ff));
ff = [O; ff];
\end{lstlisting}
The computation of the Neumann boundary condition reads
\vspace{-0.8cm}
\begin{lstlisting}
%% Assemble Neumann boundary condition
Th.elem1d = Th.bdEdge; % all boundary edges
%Th.bdEdgeIdx1 = Th.bdEdgeIdx;
%Coef = @(p) pde.Du(p)*n;
Coef = interpEdgeMat(pde.Du,Th,quadOrder);
Test = 'v.val';
ff(1:N) = ff(1:N) + assem1d(Th,Coef,Test,[],'P1',quadOrder);
\end{lstlisting}
We finally impose the Dirichlet boundary conation as
\vspace{-0.8cm}
\begin{lstlisting}
%% Apply Dirichlet boundary conditions
on = 1;
g_D = pde.g_D;
gBc = {[],g_D};
Vhvec = {'P1','P1'};
U = apply2d(on,Th,kk,ff,Vhvec,gBc); % note Vhvec
w = U(1:N);  u = U(N+1:end);
\end{lstlisting}
Note that the Dirichlet data is only for $u$, so we set \mcode{gBc\{1\} = []} in Line 4.

\begin{figure}[H]
  \centering
  \includegraphics[scale=0.5]{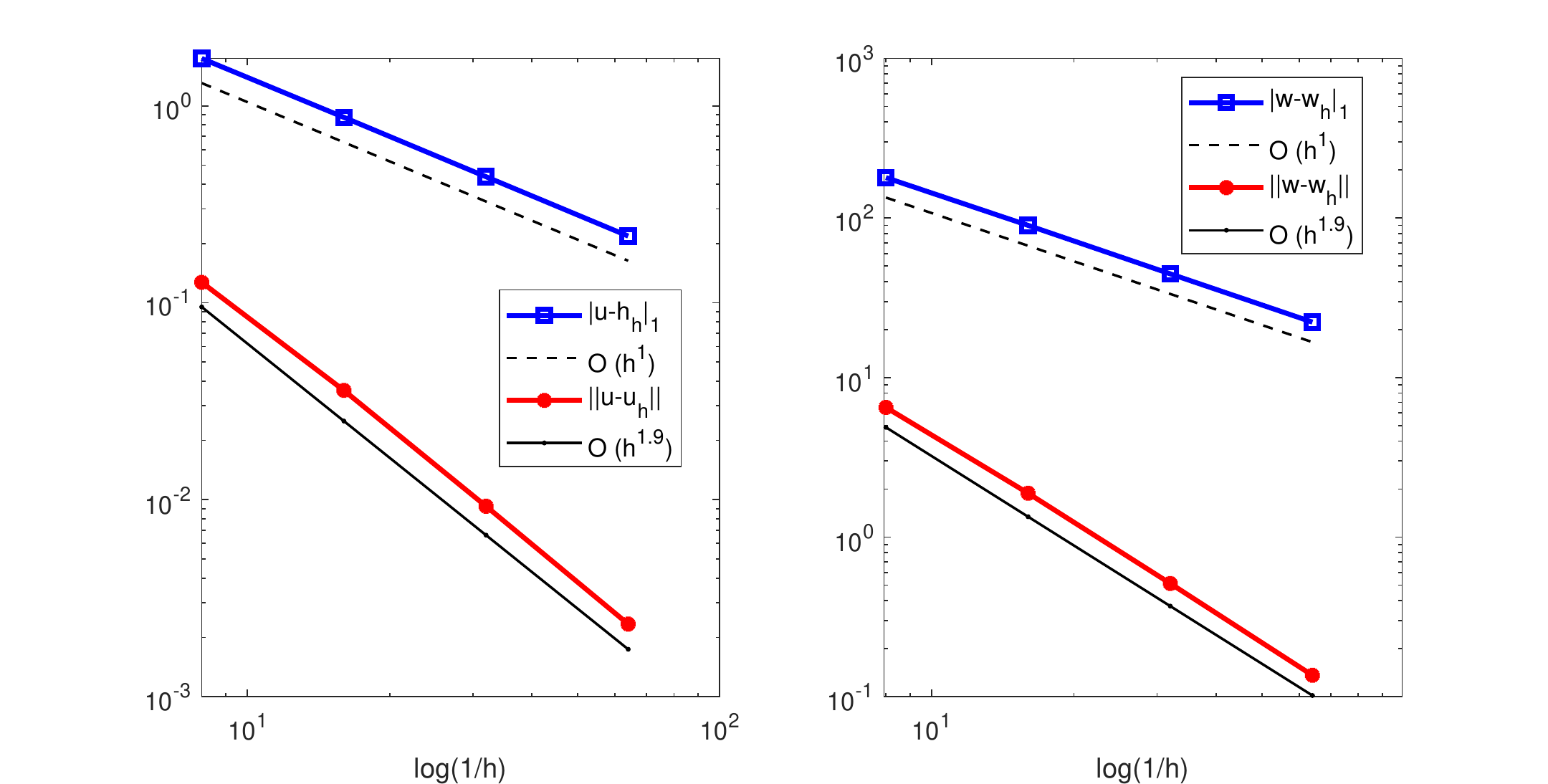}\\
  \caption{The convergence rates of the mixed FEM for the biharmonic equation ($\mathbb{P}_1$ element)}\label{fig:Biharmonic_block_rate}
\end{figure}

The convergence rates for $u$ and $w$ are shown in Fig.~\ref{fig:Biharmonic_block_rate}, from which we clearly observe the first-order and the second-order convergence in the $H^1$ norm and $L^2$ norm for both variables.

\subsubsection{The programming in vector form}

Let $u_1 = w$, $u_2 = u$, $v_1 = \phi$ and $v_2 = \psi$. Then the problem \eqref{biharmix} can be regarded as a variational problem of vector form, with $\bb u = (u_1,u_2)^T$ being the trial function and $\bb v = (v_1,v_2)^T$ being the test function. One easily finds that the mixed form \eqref{biharmix} is equivalent to the following vector form
\[\int_\Omega  ( - v_1u_1 + \nabla v_1 \cdot \nabla u_2 + \nabla v_2 \cdot \nabla u_1) \text{d}\sigma
= \int_\Omega  fv_2 \text{d}x + \int_{\partial \Omega } g_2v_1\text{d}s,\]
which is obtained by adding the two equations.

Using \mcode{int2d.m} and \mcode{int1d.m}, we can compute the vector $\bb{u}$ as follows.
\vspace{-0.8cm}
\begin{lstlisting}
%% Assemble stiffness matrix
Coef = { -1, 1, 1 };
Test  = {'v1.val', 'v1.grad', 'v2.grad'};
Trial = {'u1.val', 'u2.grad', 'u1.grad'};
kk = int2d(Th,Coef,Test,Trial,Vh,quadOrder);

%% Assemble right hand side
Coef = pde.f; Test = 'v2.val';
ff = int2d(Th,Coef,Test,[],Vh,quadOrder);

%% Assemble Neumann boundary conditions
% Get 1D mesh for boundary integrals
Th.elem1d = Th.bdEdgeType{1};
Th.elem1dIdx = Th.bdEdgeIdxType{1};

% Coef = @(p) pde.Du(p)*n;
Coef = interpEdgeMat(pde.Du,Th,quadOrder);
Test = 'v1.val';
ff = ff + int1d(Th,Coef,Test,[],Vh,quadOrder);

%% Apply Dirichlet boundary conditions
g_D = { [], pde.g_D };
on = 1;
U = apply2d(on,Th,kk,ff,Vh,g_D);
U = reshape(U,[],2);
w = U(:,1);   u = U(:,2);
\end{lstlisting}

\begin{figure}[H]
  \centering
  \subfigure[P1]{\includegraphics[scale=0.35]{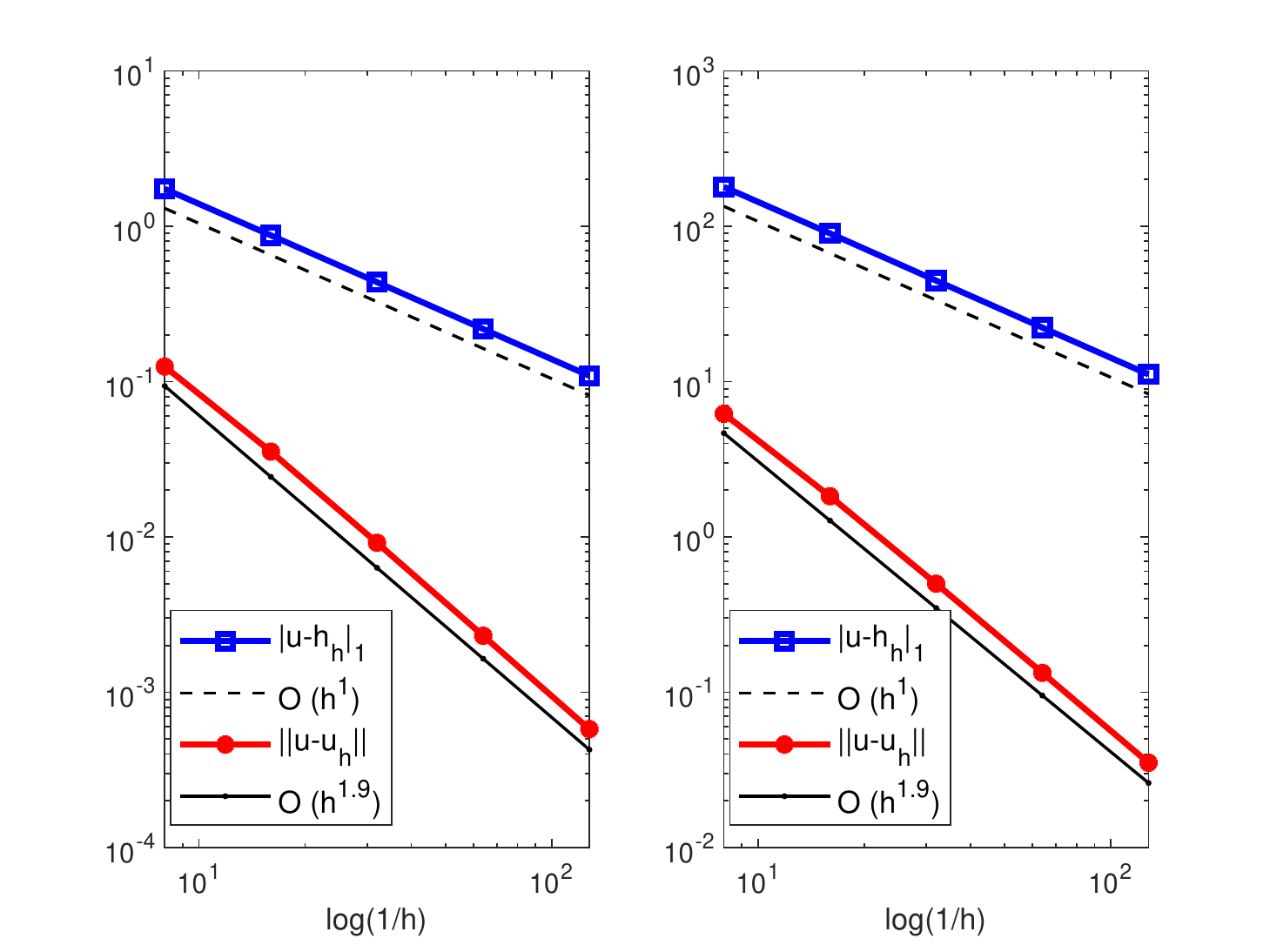}}
  \subfigure[P2]{\includegraphics[scale=0.35]{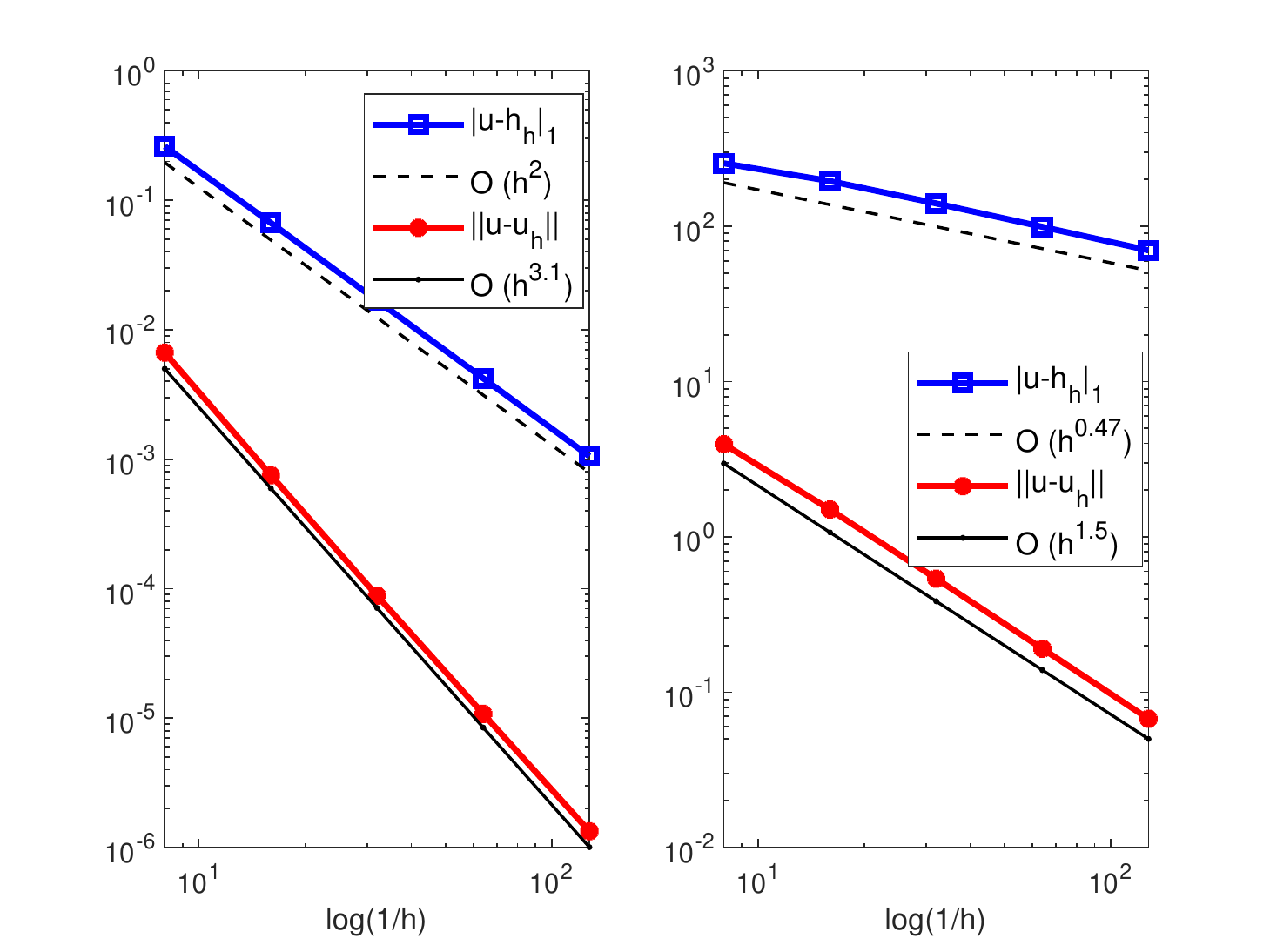}}
  \subfigure[P3]{\includegraphics[scale=0.35]{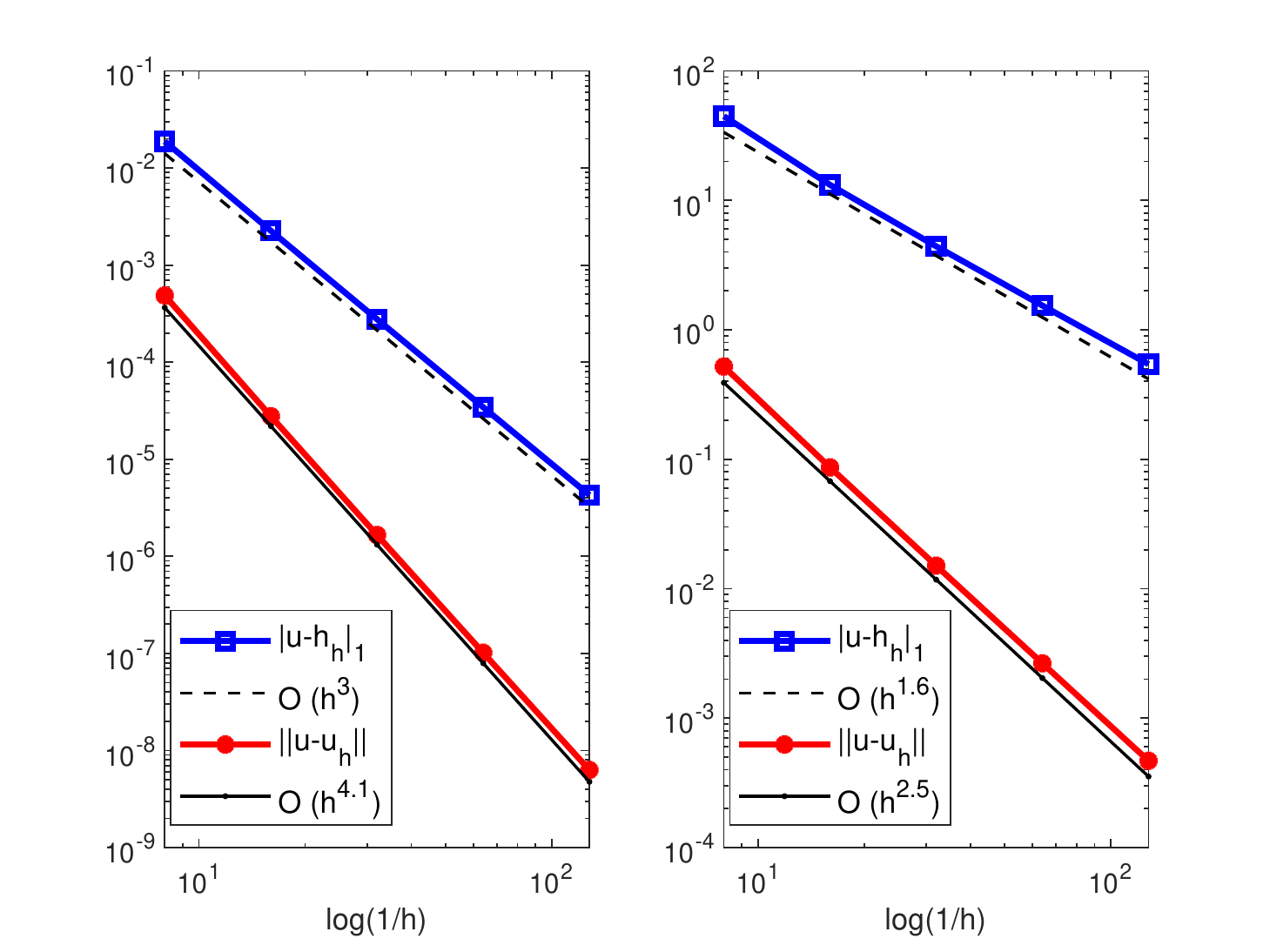}}\\
  \caption{The convergence rates for the biharmonic equation.}\label{fig:Biharmonic_rate}
\end{figure}

In the above code, \mcode{Vh} can be chosen as \mcode{\{'P1', 'P1' \}}, \mcode{\{'P2', 'P2' \}} and \mcode{\{'P3', 'P3' \}}. The results are displayed in Fig.~\ref{fig:Biharmonic_rate}. We can find that the rate of convergence for $u$ is optimal but for $w$ is sub-optimal:
\begin{itemize}
  \item[-] For linear element, optimal order for $w$ is also observed.
  \item[-] For $\mathbb{P}_2$ element, the order for $L^2$ is 1.5 and  for $H^1$ is 0.5.
  \item[-] For $\mathbb{P}_3$ element, the order for $L^2$ is 2.5 and  for $H^1$ is 1.5.
\end{itemize}
Note that our results are consistent with that given in $i$FEM. Obviously, for $\mathbb{P}_2$ and $\mathbb{P}_3$ elements, the rate of $w$ has the behaviour of $\Delta u$, which is reasonable since
$w = -\Delta u$.

\subsection{Mixed FEMs for the Stokes problem}

The Stokes problem with homogeneous Dirichlet boundary conditions is to find $(\bb{u}, p)$ such that
\[
\begin{cases}
-\nu \Delta \bb{u} - \nabla p = \bb{f} \quad & \mbox{in}~~\Omega, \\
{\rm div} \bb{u} = 0 \quad & \mbox{in}~~\Omega, \\
\bb{u} = \bb{0} \quad & \mbox{on}~~\partial \Omega.
\end{cases}
\]
Define $\bb{V} = \bb{H}_0^1(\Omega)$ and $P = L_0^2(\Omega)$. The mixed variational problem is: Find $(\bb{u}, p) \in \bb{V} \times Q$ such that
\[
\begin{cases}
a(\bb{u},\bb{v}) + b(\bb{v}, p) & = (\bb{f}, \bb{v}), \quad \bb{v}\in \bb{V}, \\
b(\bb{u},q) & = 0, \quad q\in P,
\end{cases}
\]
where
\[a(\bb{u},\bb{v}) = (\nu \nabla \bb{u}, \nabla \bb{v}), \qquad b(\bb{v}, q) = ({\rm div}\bb{v}, q).\]

Let $\mathcal{T}_h$ be a shape regular triangulation  of $\Omega $. We consider the conforming finite element discretizations: $\bb{V}_h \subset \bb{V}$ and $P_h \subset P$. Typical pairs $(\bb{V}_h,P_h)$ of stable finite element spaces include: MINI element, Girault-Raviart element and $\mathbb{P}_k-\mathbb{P}_{k - 1}$ elements. For the last one, a special example is the $\mathbb{P}_2-\mathbb{P}_1$ element, also known as the Taylor-Hood element, which is the one under consideration.

The FEM is to find $(\bb u_h,p_h) \in \bb{V}_h \times P_h$ such that
\begin{equation}\label{Stokefem}
\begin{cases}
  a(\bb u_h,\bb v) + b(\bb v,p_h) = F(\bb v), \quad & \bb v \in \bb{V}_h,  \\
  b(\bb u_h,q) = 0,\quad & q \in P_h.
\end{cases}
\end{equation}
The problem \eqref{Stokefem} can be solved either by discretizing it directly into a system of equations (a saddle point problem), or by adding the two equations, as done for the biharmonic equation. We consider the latter one: Find $(\bb u_h,p_h) \in \bb{V}_h \times L^2(\Omega)$ such that
\begin{equation}\label{Stokefem1}
a(\bb u_h,\bb v) + b(\bb v,p_h) + b(\bb u_h,q) - \varepsilon (p_h,q) = F(\bb v), \bb v \in \bb{V}_h, \quad  \in L^2(\Omega ),
\end{equation}
where $\varepsilon$ is a small parameter to ensure stability.

The bilinear form can be assembled as follows.
\vspace{-0.8cm}
\begin{lstlisting}
Vh = {'P2','P2','P1'}; quadOrder = 5;

%% Assemble stiffness matrix
vstr = {'v1','v2','q'}; ustr = {'u1','u2','p'};
% [v1,v2,q] = [v1,v2,v3], [u1,u2,p] = [u1,u2,u3]
%   a1(v1,u1) + a2(v2,u2)
% + b1(v1,p)  + b2(v2,p)
% + b1(q,u1)  + b2(q,u2) - eps*(q,p)
eps = 1e-10;
Coef = { 1, 1,  -1,-1,  -1,-1,  -eps};
Test  = {'v1.grad', 'v2.grad', 'v1.dx',  'v2.dy',  'q.val', 'q.val', 'q.val'};
Trial = {'u1.grad', 'u2.grad', 'p.val',  'p.val', 'u1.dx',  'u2.dy',  'p.val'};
[Test,Trial] = getStdvarForm(vstr, Test,  ustr, Trial); % [u1,u2,p] --> [u1,u2,u3]
[kk,info] = int2d(Th,Coef,Test,Trial,Vh,quadOrder);
\end{lstlisting}
Note that the symbols \mcode{p} and \mcode{q} correspond to \mcode{u3} and \mcode{v3}, which is realized by using the subroutine \mcode{getStdvarForm.m}.

The computation of the right-hand side reads
\vspace{-0.8cm}
\begin{lstlisting}
%% Assemble right hand side
trf = eye(2);
Coef1 = @(pz) pde.f(pz)*trf(:, 1);  Coef2 = @(pz) pde.f(pz)*trf(:, 2);
Coef = {Coef1, Coef2};
Test = {'v1.val', 'v2.val'};
ff = int2d(Th,Coef,Test,[],Vh,quadOrder);
\end{lstlisting}
In this case, one can not use \mcode{Coef = pde.f} and \mcode{Test = 'v.val'} instead since $\bb{v} = [v_1,v_2,v_3]^T$ has three components.

We impose the Dirichlet boundary conditions as follows.
\vspace{-0.8cm}
\begin{lstlisting}
%% Apply Dirichlet boundary conditions
tru = eye(2);
g_D1 = @(pz) pde.g_D(pz)*tru(:, 1);
g_D2 = @(pz) pde.g_D(pz)*tru(:, 2);
g_D = {g_D1, g_D2, []};
on = 1;
U = apply2d(on,Th,kk,ff,Vh,g_D);
NNdofu = info.NNdofu;
id1 =  NNdofu(1);  id2 = NNdofu(1)+ NNdofu(2);
uh = [ U(1:id1), U(id1+1:id2) ]; % uh = [u1h, u2h]
ph = U(id2+1:end);
\end{lstlisting}
Note that \mcode{g_D\{3\} = []} since no constraints are imposed on $p$.

\begin{example}
Let $\Omega = (0,1)^2$. We choose the load term $\bb{f}$ in such a way that the analytical solution is
\[\bb{u}(x,y) =
\begin{bmatrix}
-2^8(x^2-2x^3+x^4)(2y-6y^2+4y^3)\\
2^8(2x-6x^2+4x^3)(y^2-2y^3+y^4)
 \end{bmatrix}\]
 and $p(x,y) =  -2^8(2-12x+12x^2)(y^2-2y^3+y^4)$.
\end{example}

\begin{figure}[H]
  \centering
  \includegraphics[scale=0.45]{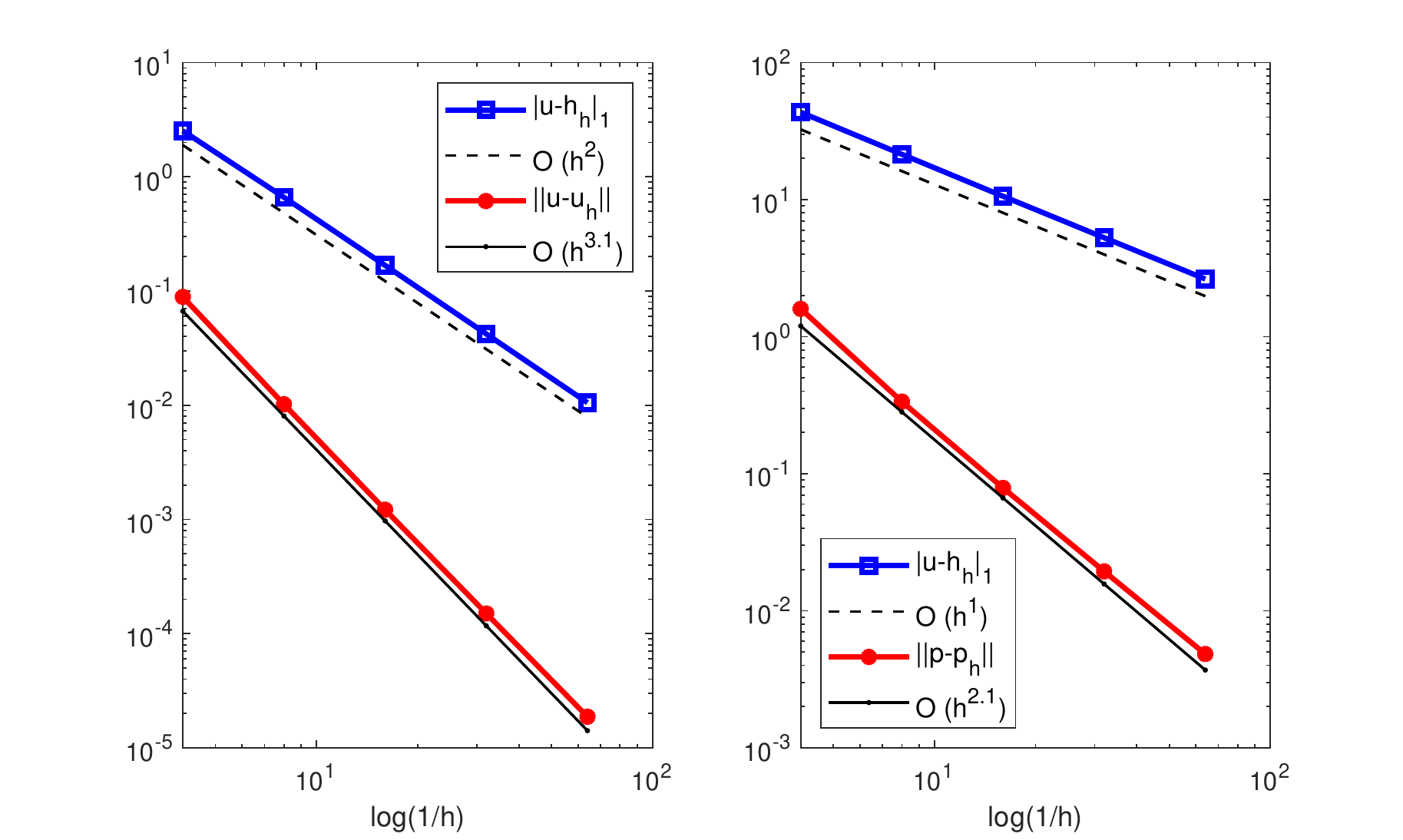}\\
  \caption{Convergence rates of the Taylor-Hood element for the Stokes problem}\label{fig:Stokes_rate}
\end{figure}

The results are displayed in Fig.~\ref{fig:Stokes_rate} and Tab.~\ref{tab:Stokes_rate}, from which we observe the optimal rates of convergence both for $u$ and $p$.

\begin{table}[!htb]
 \centering
 \caption{The discrete errors for the Stokes problem} \label{tab:Stokes_rate}
\begin{tabular}{rcccccccccc}
  \hline
  $\sharp$Dof   &    $h$    &    $\|u-u_h\|$ &   $|u-u_h|_1$  &  $ \|p-p_h\|$ \\
  \hline
      32  &  2.500e-01  &  8.88464e-02  &  2.52940e+00   & 1.59802e+00 \\
     128  &  1.250e-01  &  1.01868e-02  &  6.62003e-01   & 3.36224e-01 \\
     512  &  6.250e-02  &  1.21537e-03  &  1.67792e-01   & 7.88512e-02 \\
    2048  &  3.125e-02  &  1.50235e-04  &  4.21077e-02   & 1.94079e-02 \\
    8192  &  1.562e-02  &  1.87368e-05  &  1.05374e-02   & 4.83415e-03 \\
  \hline
\end{tabular}
\end{table}

\subsection{Time-dependent problems}

As an example, we consider the heat equation:
\[\begin{cases}
u_t - \Delta u = f  & \mbox{in} ~~\Omega, \\
u(x,y,0) = u_0(x,y) & \mbox{in} ~~\Omega, \\
u = g_D  & \mbox{on} ~~\Gamma_D,  \\
\partial_nu =g_N &  \mbox{on} ~~\Gamma_N.
\end{cases}\]
After applying the backward Euler discretization in time, we shall seek $u^n(x,y)$ satisfying for all $v\in H_0^1(\Omega)$:
\[\int_\Omega \Big( \frac{u^n-u^{n-1}}{\Delta t} v + \nabla u^n \cdot \nabla v \Big) {\rm d}\sigma
=  \int_\Omega f^n v {\rm d}\sigma + \int_{\Gamma_N} g_N^n v {\rm d}s.\]

We fist generate a mesh and compute the mesh information.
\vspace{-0.8cm}
\begin{lstlisting}
%% Mesh
% generate mesh
Nx = 10;
[node,elem] = squaremesh([0 1 0 1],1/Nx);
% mesh info
bdStr = 'x==0'; % Neumann
Th = FeMesh2d(node,elem,bdStr);
% time
Nt = Nx^2;
t = linspace(0,1,Nt+1)';  dt = t(2)-t(1);
\end{lstlisting}
The PDE data is given by
\vspace{-0.8cm}
\begin{lstlisting}
%% PDE data
pde = heatData();
\end{lstlisting}
The exact solution is chosen as $u = \sin(\pi x)\sin(y)\e^{-t}$.

For fixed $\Delta t$, the bilinear form gives the same stiffness matrix in each iteration.
\vspace{-0.8cm}
\begin{lstlisting}
%% Backward Euler
u0 = @(p) pde.uexact(p,t(1));
uh0 = interp2d(u0,Th,Vh); % dof vector
uf = zeros(Nt+1,2);  % record solutions at p-th point
p = 2*Nx;
uf(1,:) = [uh0(p),uh0(p)];
for n = 1:Nt
    % Linear form
    fun = @(p) pde.f(p, t(n+1));
    Coef = fun;  Test = 'v.val';
    ff = assem2d(Th,Coef,Test,[],Vh,quadOrder);
    Coef = uh0/dt;
    ff = ff + assem2d(Th,Coef,Test,[],Vh,quadOrder);

    % Neumann boundary condition
    if ~isempty(bdStr)
        Th.elem1d = Th.bdEdgeType{1};
        Th.elem1dIdx = Th.bdEdgeIdxType{1};
        fun = @(p) pde.Du(p,t(n+1));
        Coef = interpEdgeMat(fun,Th,quadOrder);
        ff = ff + assem1d(Th,Coef,Test,[],Vh,quadOrder);
    end

    % Dirichlet boundary condition
    ue = @(p) pde.uexact(p,t(n+1));
    on = 2 - 1*isempty(bdStr); % on = 1 for the whole boundary if bdStr = []
    uh = apply2d(on,Th,kk,ff,Vh,ue);

    % Record
    uhe = interp2d(ue,Th,Vh);
    uf(n+1,:) = [uhe(p), uh(p)];

    % Update
    uh0 = uh;
end
\end{lstlisting}

\begin{figure}[H]
  \centering
  \includegraphics[scale=0.5]{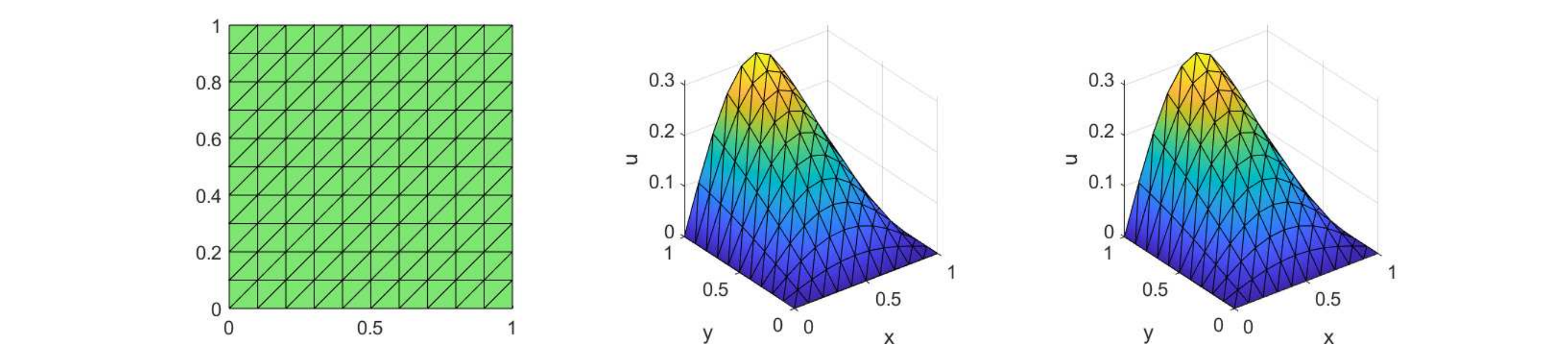}\\
  \caption{Exact and numerical solutions for the heat equation}\label{fig:heat_val}
\end{figure}

\begin{figure}[H]
  \centering
  \includegraphics[scale=0.4]{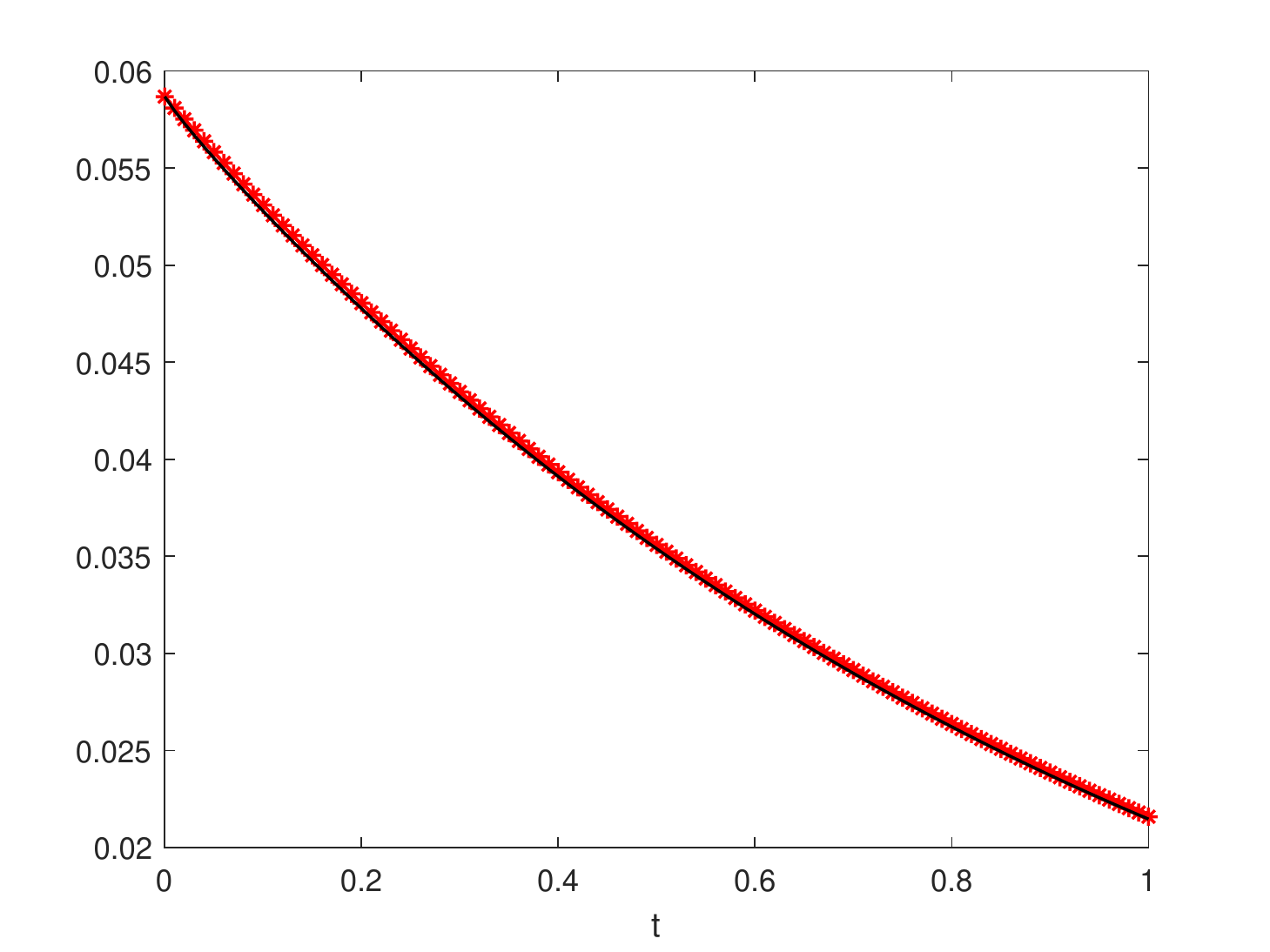}\\
  \caption{Exact and numerical solutions at a fixed point for the heat equation}\label{fig:heat_valx}
\end{figure}

In the above code, we record solutions at $p$-th point, where $p= 2N_x$. The results are shown in Fig.~\ref{fig:heat_val} and Fig.~\ref{fig:heat_valx}.

\begin{figure}[H]
  \centering
  \subfigure[P1]{\includegraphics[scale=0.35]{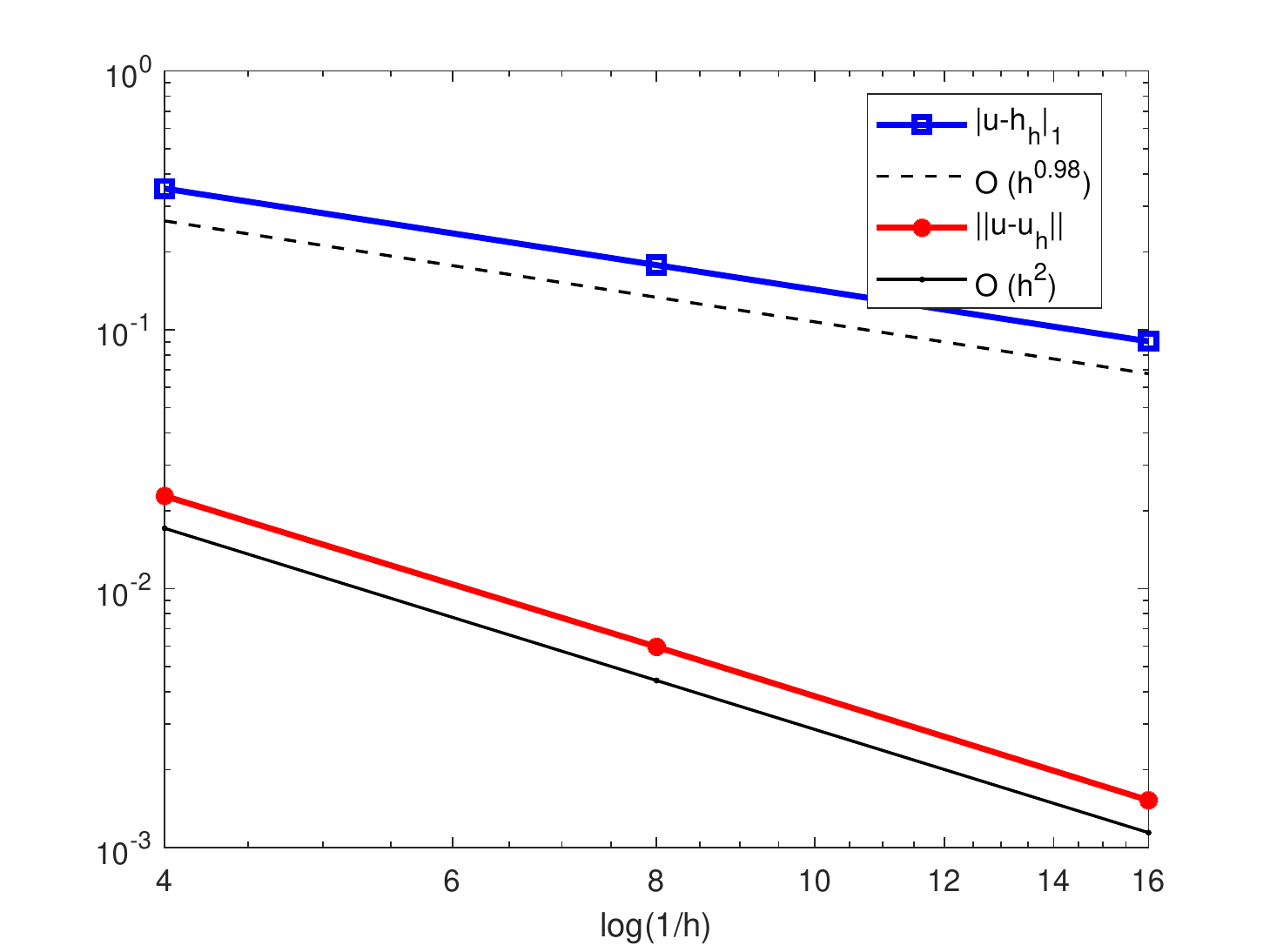}}
  \subfigure[P2]{\includegraphics[scale=0.35]{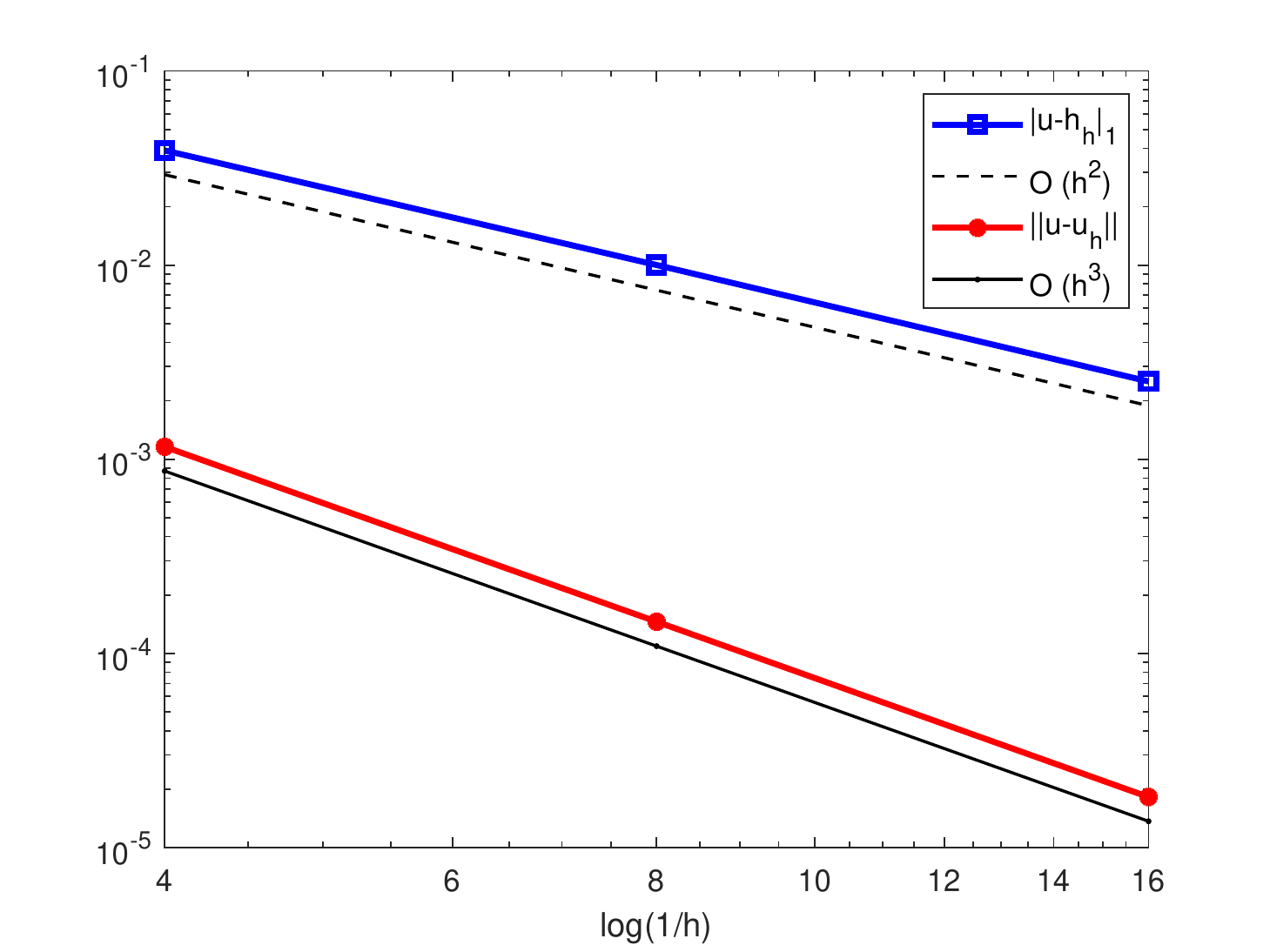}}
  \subfigure[P3]{\includegraphics[scale=0.35]{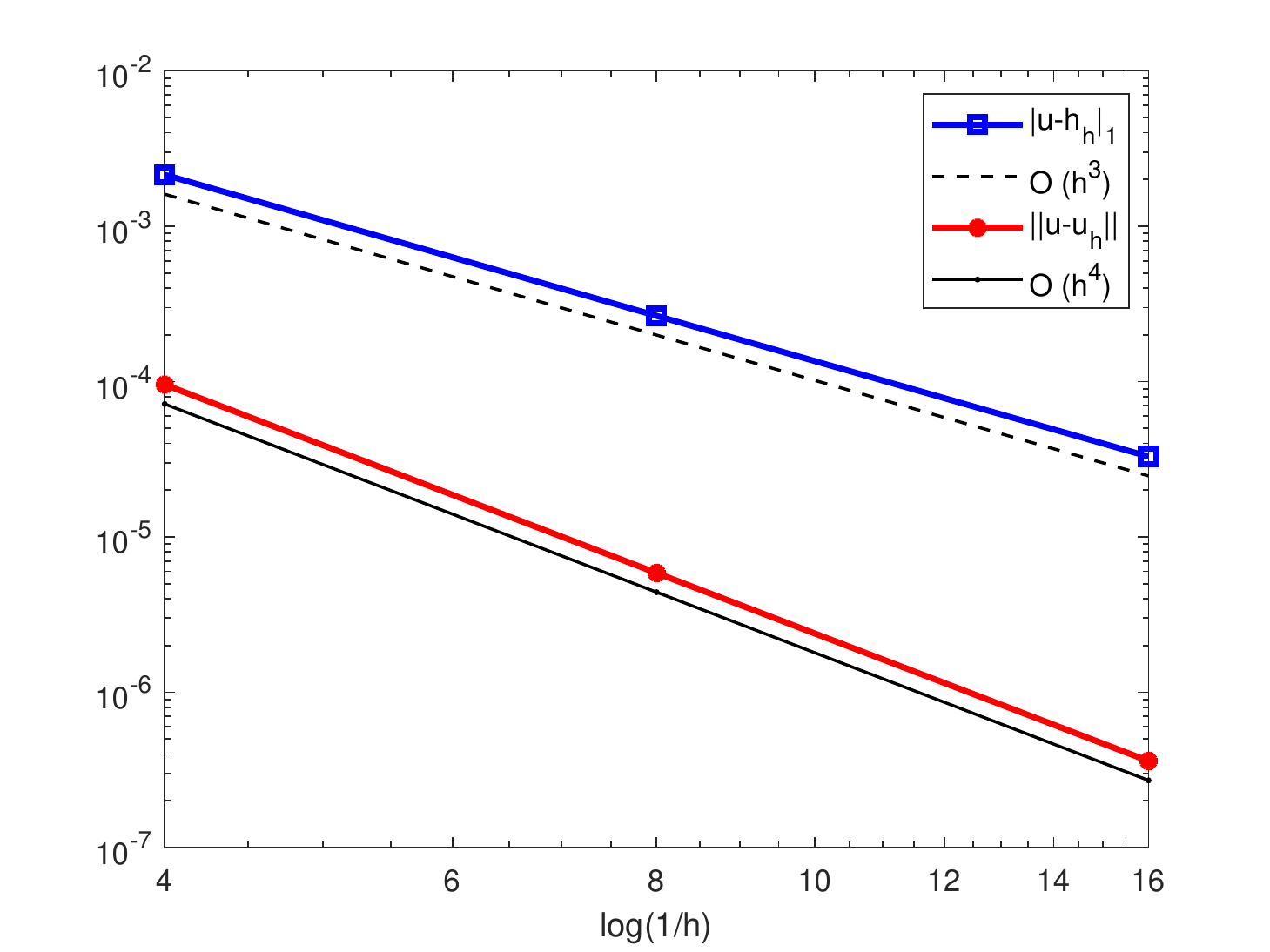}}\\
  \caption{Convergence rates for the heat equation.}\label{fig:heat_rate}
\end{figure}

It is well-known that the error behaves as $\mathcal{O}(\Delta t + h^{k+1})$ for the $\mathbb{P}_k$-Lagrange element. To test the convergence order w.r.t $h$, one can set $\Delta t = \mathcal{O}(h^{k+1})$. The numerical results are consistent with the theoretical prediction as shown in Fig.~\ref{fig:heat_rate}.

\section{Examples in FreeFEM Documentation} \label{sect:freefem}

We in this section present five examples given in FreeFEM. Many more examples can be found or will be added in the example folder in varFEM.

\subsection{Membrane}

This is an exmple given in FreeFEM Documentation: Release 4.6 (see Subsection 2.3 - Membrane).

The equation is simply the Laplace equation, where the region is an ellipse with the length of the semimajor axis $a = 2$, and unitary the semiminor axis. The mesh on such a domain can be generated by using the pdetool:
\vspace{-0.8cm}
\begin{lstlisting}
%% Mesh
% ellipse with a = 2, b = 1
a = 2; b = 1;
g = ellipseg(a,b);
[p,e,t] = initmesh(g,'hmax',0.2);
[p,e,t] = refinemesh(g,p,e,t);
node = p'; elem = t(1:3,:)';
figure,
subplot(1,2,1),
showmesh(node,elem);
% bdStr
bdNeumann = 'y<0 & x>-sin(pi/3)'; % string for Neumann
% mesh info
Th = FeMesh2d(node,elem,bdNeumann);
\end{lstlisting}

The Neumann boundary condition is imposed on
\[\Gamma_2 = \{ (x,y): x = a \cos t,  ~~y = b \sin t, \quad t \in (\theta, 2\pi) \}, \qquad \theta = 4\pi/3,\]
which can be identified by setting \mcode{bdStr = 'y<0 & x>-sin(pi/3)'} as done in the code.

\begin{figure}[H]
  \centering
  \subfigure[Nodal values]{\includegraphics[scale=0.45,trim=0 40 0 40,clip]{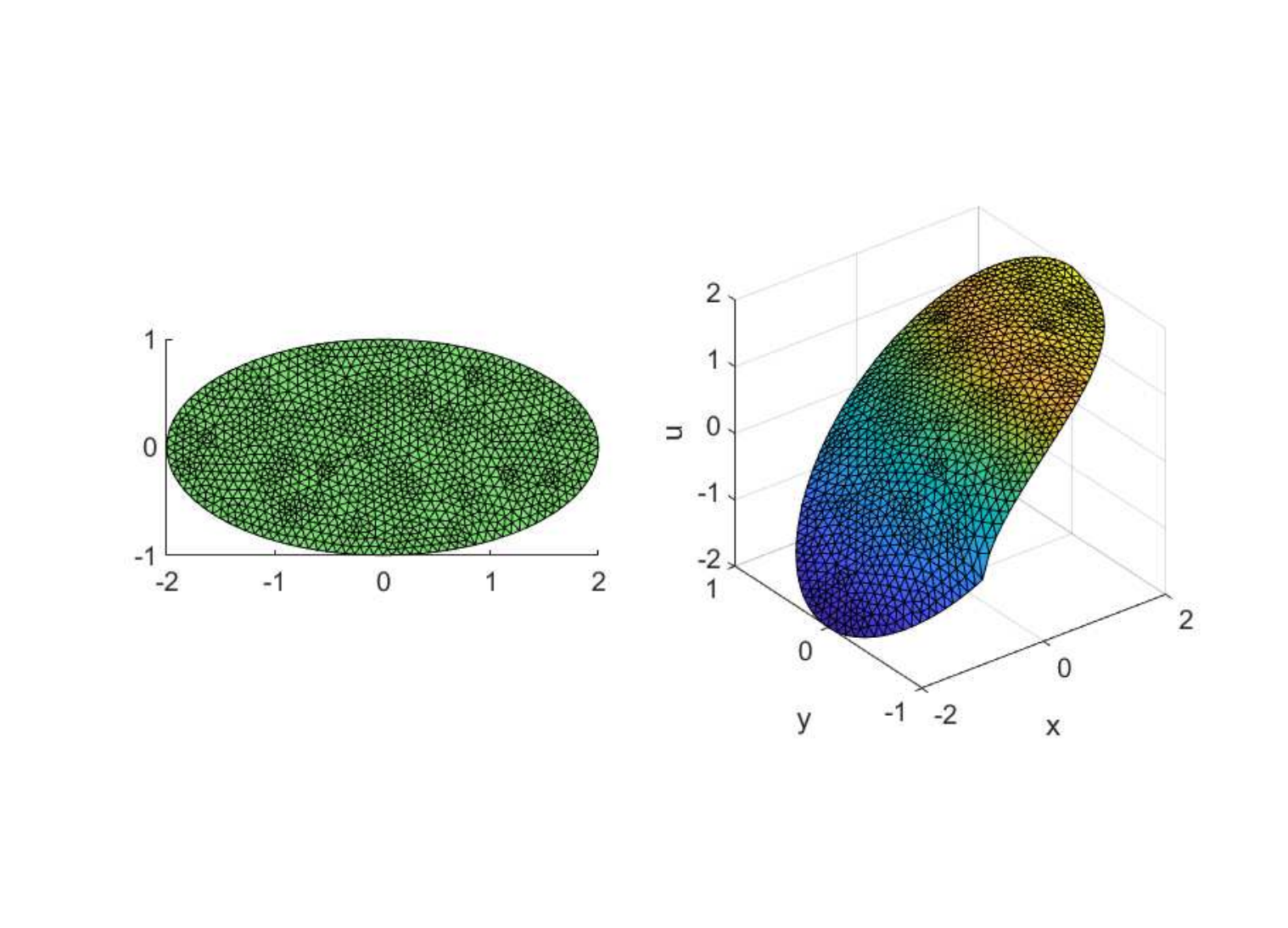}}
  \subfigure[Level lines]{\includegraphics[scale=0.45,trim=0 40 0 40,clip]{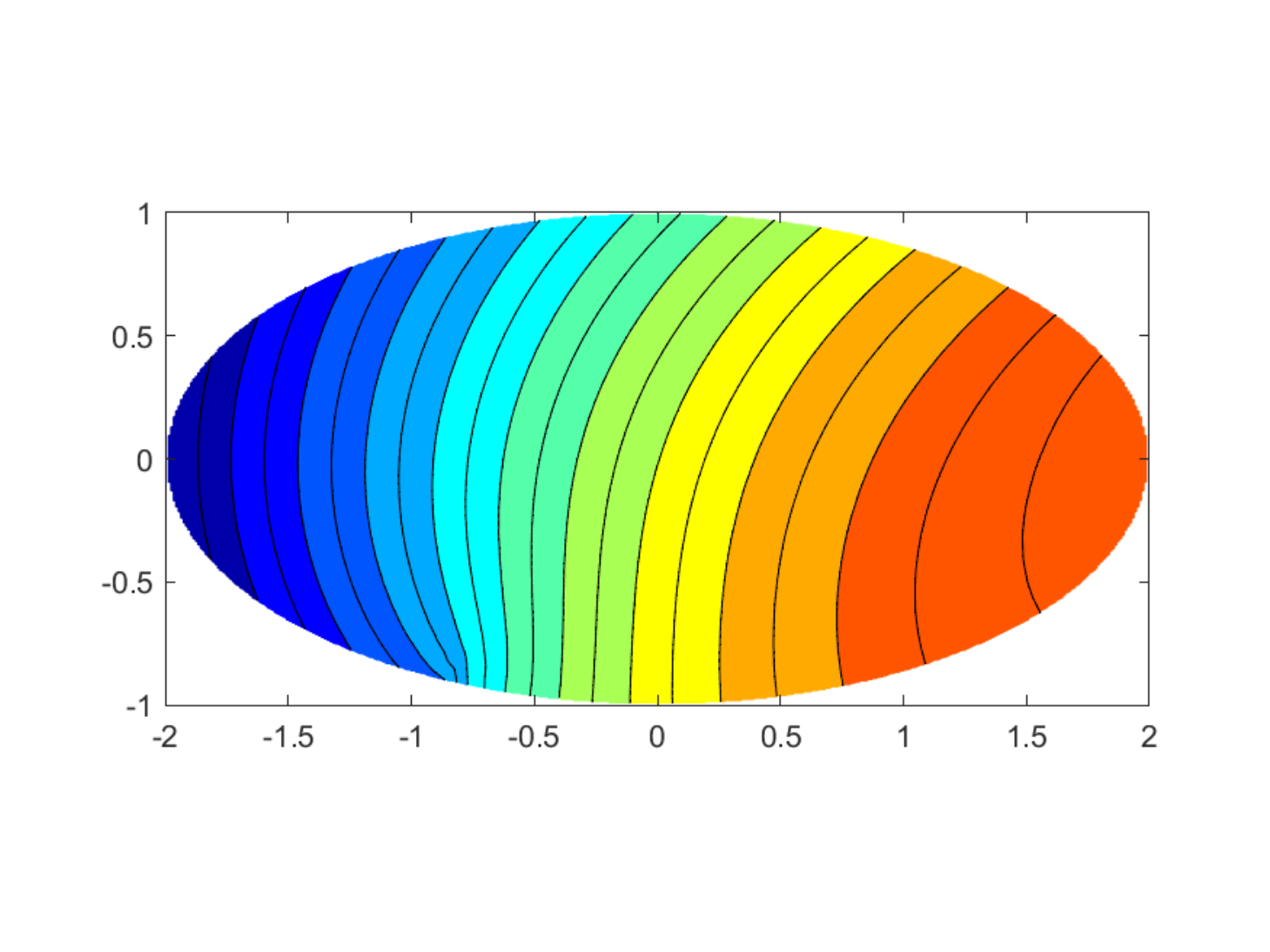}}
  \\
  \caption{Membrane deformation}\label{fig:membrane_val}
\end{figure}

The remaining implementation is very simple, so we omit the details. The nodal values are shown in Fig.~\ref{fig:membrane_val}a.

We can also plot the level lines of the membrane deformation, as given in Fig.~\ref{fig:membrane_val}b. Note that the contour figure is obtained by interpolating the finite element function to a two-dimensional cartesian grid (within the mesh). The interpolated values can be created by using the Matlab build-in function \mcode{pdeInterpolant.m} in the pdetool. However, the build-in function seems not efficient. For this reason, we provide a new realization of the interpolant, named \mcode{varInterpolant2d.m}. With this function, we give a subroutine \mcode{varcontourf.m} to draw a contour plot.

\subsection{Heat exchanger}

\begin{figure}[H]
  \centering
  \includegraphics[scale=0.45]{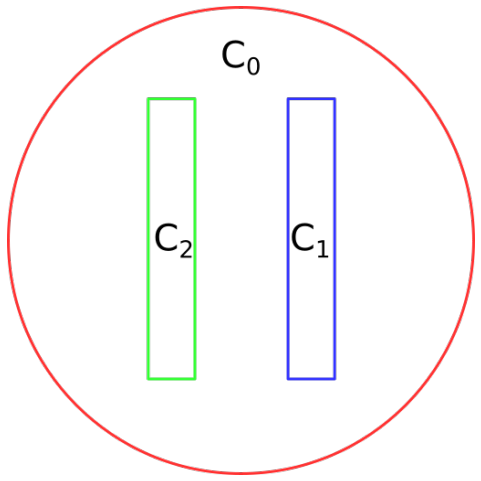}\\
  \caption{Heat exchanger geometry (see \cite{FreeFEM})}\label{fig:heatexGeo}
\end{figure}

The geometry is shown in Fig.~\ref{fig:heatexGeo}, where $C_1$ and $C_2$ are two thermal conductors within an enclosure $C_0$. The temperature $u$ satisfies
\[\nabla \cdot ( \kappa \nabla u) = 0 \quad \mbox{in}~~\Omega.\]
\begin{itemize}
  \item The first conductor is held at a constant temperature $u_1 = 100^\circ {\rm C}$, and the border of enclosure $C_0$ is held at temperature $20^\circ {\rm C}$. This means the domain is $\Omega = C_0 \backslash C_1$, and the boundaries consist of $\partial C_0$ and $\partial C_1$.
  \item The conductor $C_2$ has a different thermal conductivity than the enclosure $C_0$: $\kappa_0 = 1$ and $k_2 = 3$.
\end{itemize}

We use the mesh generated by FreeFEM, which is saved in a .msh file named \mcode{meshdata\_heatex.msh}. The command in FreeFEM can be written as
\vspace{-0.8cm}
\begin{lstlisting}
savemesh(Th, "meshdata_heatex.msh");
\end{lstlisting}
We read the basic data structures \mcode{node} and \mcode{elem} via a self-written function:
\vspace{-0.8cm}
\begin{lstlisting}
[node,elem] = getMeshFreeFEM('meshdata_heatex.msh');
\end{lstlisting}
Then the mesh data can be computed as
\vspace{-0.8cm}
\begin{lstlisting}
% bdStr
C0 = 'x.^2 + y.^2 > 3.8^2'; % 1
bdStr = C0;
% mesh info
Th = FeMesh2d(node,elem,bdStr);
\end{lstlisting}
Note that the remaining boundary is $\partial C_1$, hence the boundaries of $C_0$ and the first conductor are labelled as 1 and 2, respectively.

The coefficient $\kappa$ can be written as
\vspace{-0.8cm}
\begin{lstlisting}
%% PDE data
kappa = @(p) 1 + 2*(p(:,1)<-1).*(p(:,1)>-2).*(p(:,2)<3).*(p(:,2)>-3);  % p = [x,y]
\end{lstlisting}
And the bilinear form and the linear form are assembled as follows.
\vspace{-0.8cm}
\begin{lstlisting}
%% Bilinear form
Vh = 'P1'; quadOrder = 5;
Coef = kappa;
Test = {'v.grad'};  Trial = {'u.grad'};
kk = assem2d(Th,Coef,Test,Trial,Vh,quadOrder);

%% Linear form
ff = zeros(size(kk,1),1);
\end{lstlisting}

The Dirichlet boundary conditions are imposed in the following way:
\vspace{-0.8cm}
\begin{lstlisting}
%% Dirichlet boundary conditions
on = [1,2];
gBc1 = @(p) 20+0*p(:,1);
gBc2 = @(p) 100+0*p(:,1);
uh = apply2d(on,Th,kk,ff,Vh,gBc1,gBc2);
\end{lstlisting}
Here, \mcode{gBc1} is for $\partial C_0$ or \mcode{on(1)}, and \mcode{gBc2} is for $\partial C_1$ or \mcode{on(2)}.

In FreeFEM, the numerical solution can be outputted using the command \mcode{ofstream}:
\vspace{-0.8cm}
\begin{lstlisting}
ofstream file("sol_heatex.txt");
file<<u[]<<endl;
\end{lstlisting}
The saved information is then loaded by
\vspace{-0.8cm}
\begin{lstlisting}
uff = solFreeFEM('sol_heatex.txt');
\end{lstlisting}

\begin{figure}[H]
  \centering
  \subfigure[varFEM]{\includegraphics[scale=0.35]{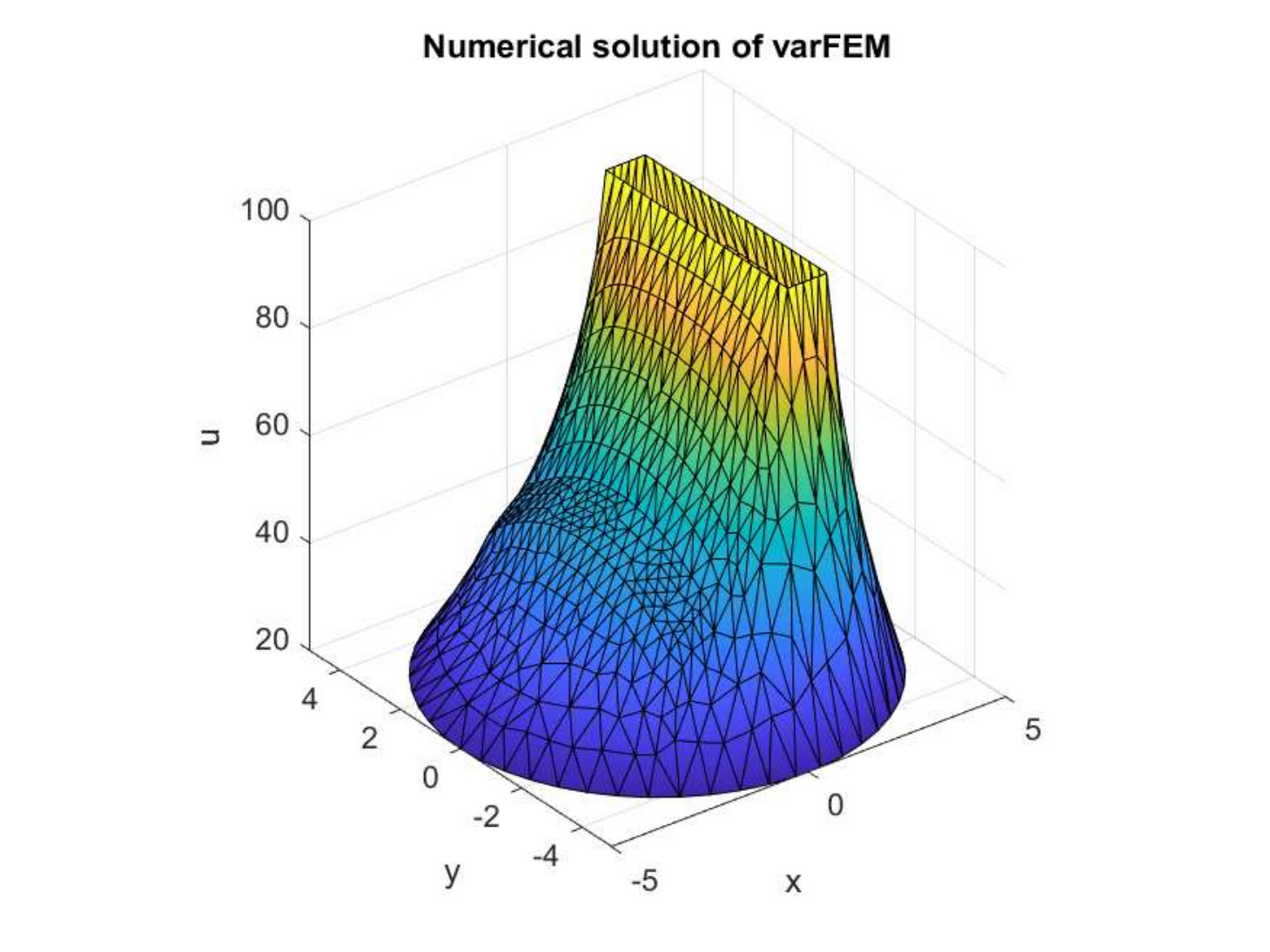}}
  \subfigure[FreeFEM]{\includegraphics[scale=0.35]{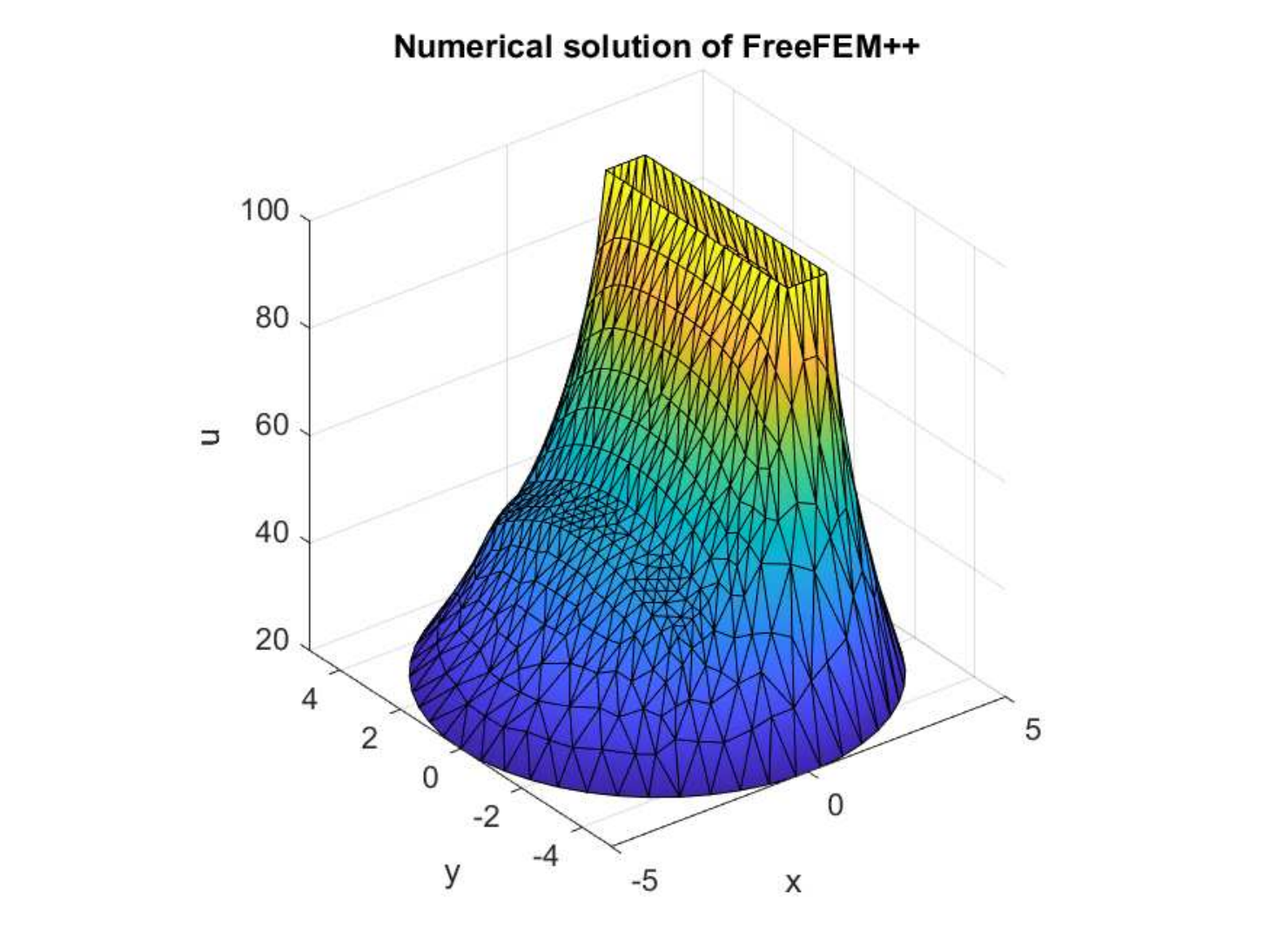}}
  \subfigure[Level lines]{\includegraphics[scale=0.35]{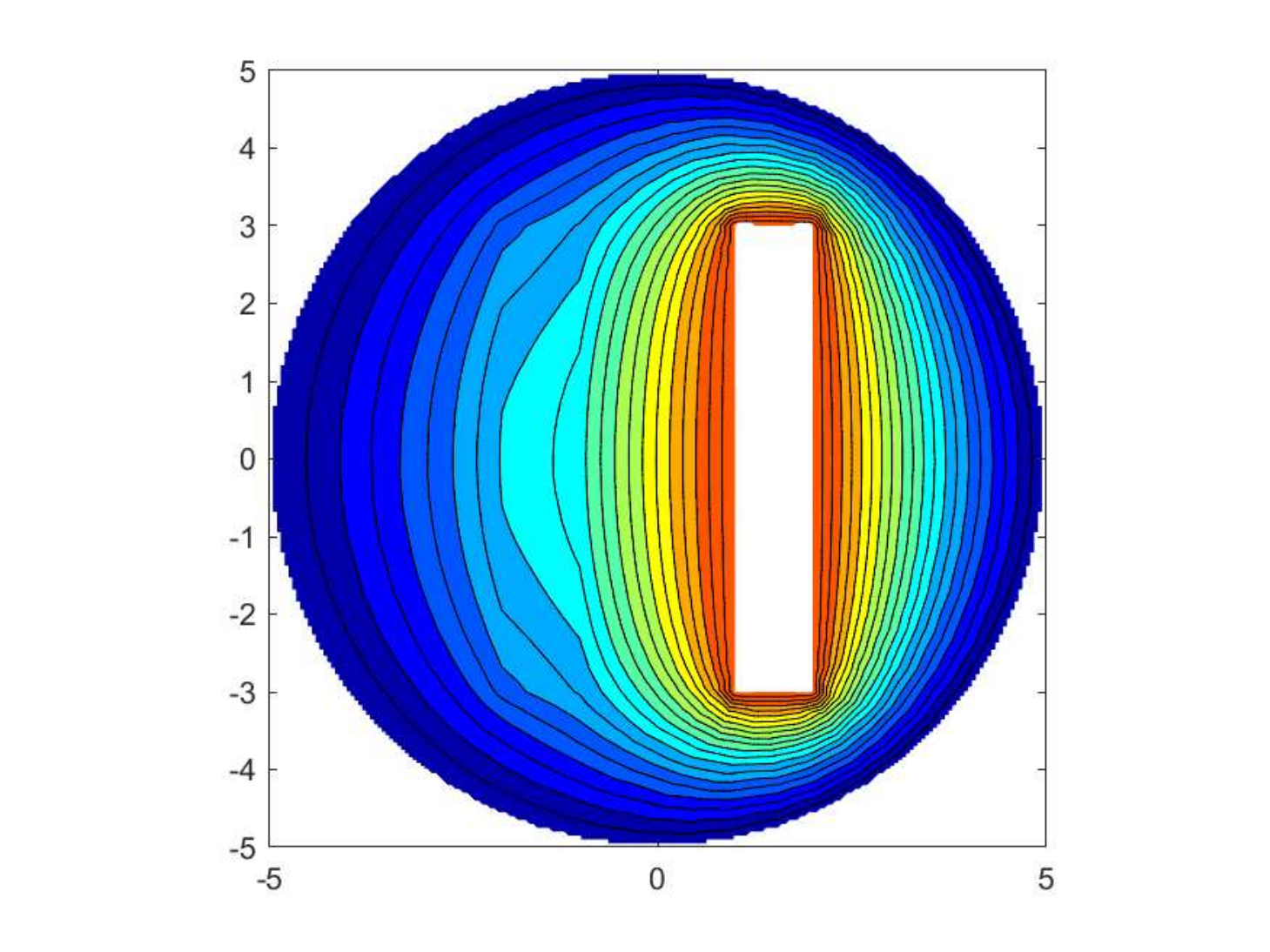}}\\
  \caption{Numerical solutions given by varFEM and FreeFEM}\label{fig:heatexval}
\end{figure}

The results are shown in Fig.~\ref{fig:heatexval}, from which we observe that the varFEM solution is well matched with the one given by FreeFEM.

\subsection{Airfoil}

This is an exmple given in FreeFem Documentation: Release 4.6 (Subsection 2.7 - Irrotational Fan Blade Flow and Thermal effects).

\begin{figure}[H]
  \centering
   \includegraphics[scale=0.4,trim=0 60 0 60,clip]{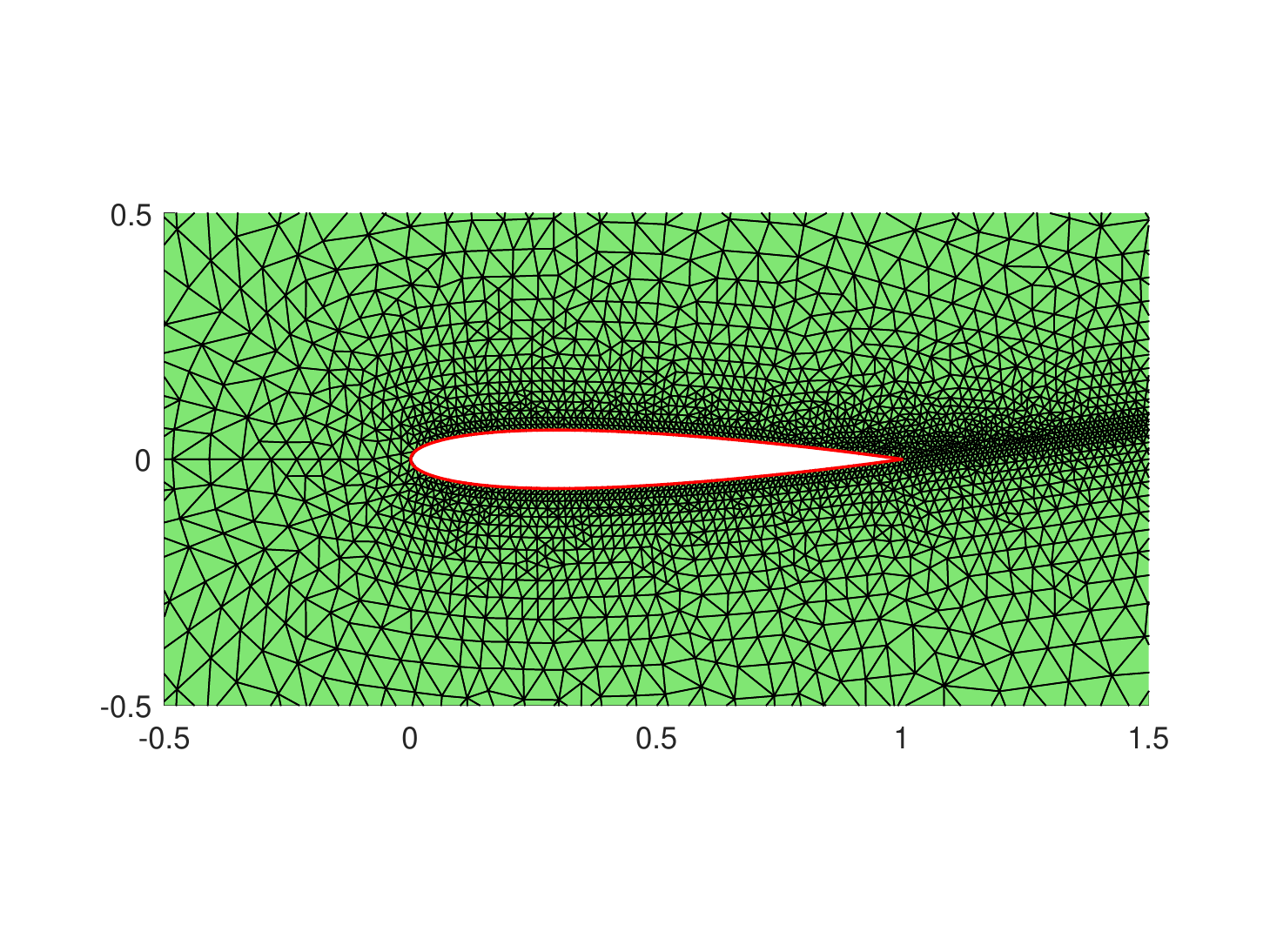}
   \caption{Mesh zoomed around the NACA0012 airfoil}\label{fig:naca0012}
\end{figure}

 Consider a wing profile $S$ (the NACA0012 Airfoil) in a uniform flow. Infinity will be represented by a large circle $C$ where the flow is assumed to be of uniform velocity. The domain is outside $S$, with the mesh shown in Fig.~\ref{fig:naca0012}. The NACA0012 airfoil is a classical wing profile in aerodynamics, whose equation for the upper surface is
 \[y = 0.17735\sqrt{x} - 0.075597x - 0.212836x^2 + 0.17363x^3 - 0.06254x^4.\]
With this equation, we can generate a mesh using the Matlab pdetool, as included in varFEM. The function is \mcode{mesh\_naca0012.m}. For comparison, we use the mesh generated by FreeFEM.

The programming is very simple, given by
\vspace{-0.8cm}
\begin{lstlisting}
%% Parameters
theta = 8*pi/180;
lift = theta*0.151952/0.0872665; % lift approximation formula
uinfty1 = cos(theta);   uinfty2 = sin(theta);

%% Mesh
[node,elem] = getMeshFreeFEM('meshdata_airfoil.msh');
% mesh info
bdStr = 'x.^2 + y.^2 > 4.5^2'; % 1-C
Th = FeMesh2d(node,elem,bdStr);

%% Bilinear form
Vh = 'P2';  quadOrder = 7;
Coef  = 1;
Test  = 'v.grad';
Trial = 'u.grad';
kk = assem2d(Th,Coef,Test,Trial,Vh,quadOrder);

%% Linear form
ff = zeros(size(kk,1),1);

%% Dirichlet boundary conditions
on = [1,2];
gBc1 = @(p) uinfty1*p(:,2) - uinfty2*p(:,1); % on 1-C
gBc2 = @(p) -lift + 0*p(:,1);  % on 2-S
uh = apply2d(on,Th,kk,ff,Vh,gBc1,gBc2);
\end{lstlisting}

\begin{figure}[H]
  \centering
  \subfigure[varFEM solution]{\includegraphics[scale=0.4,trim=0 60 0 60,clip]{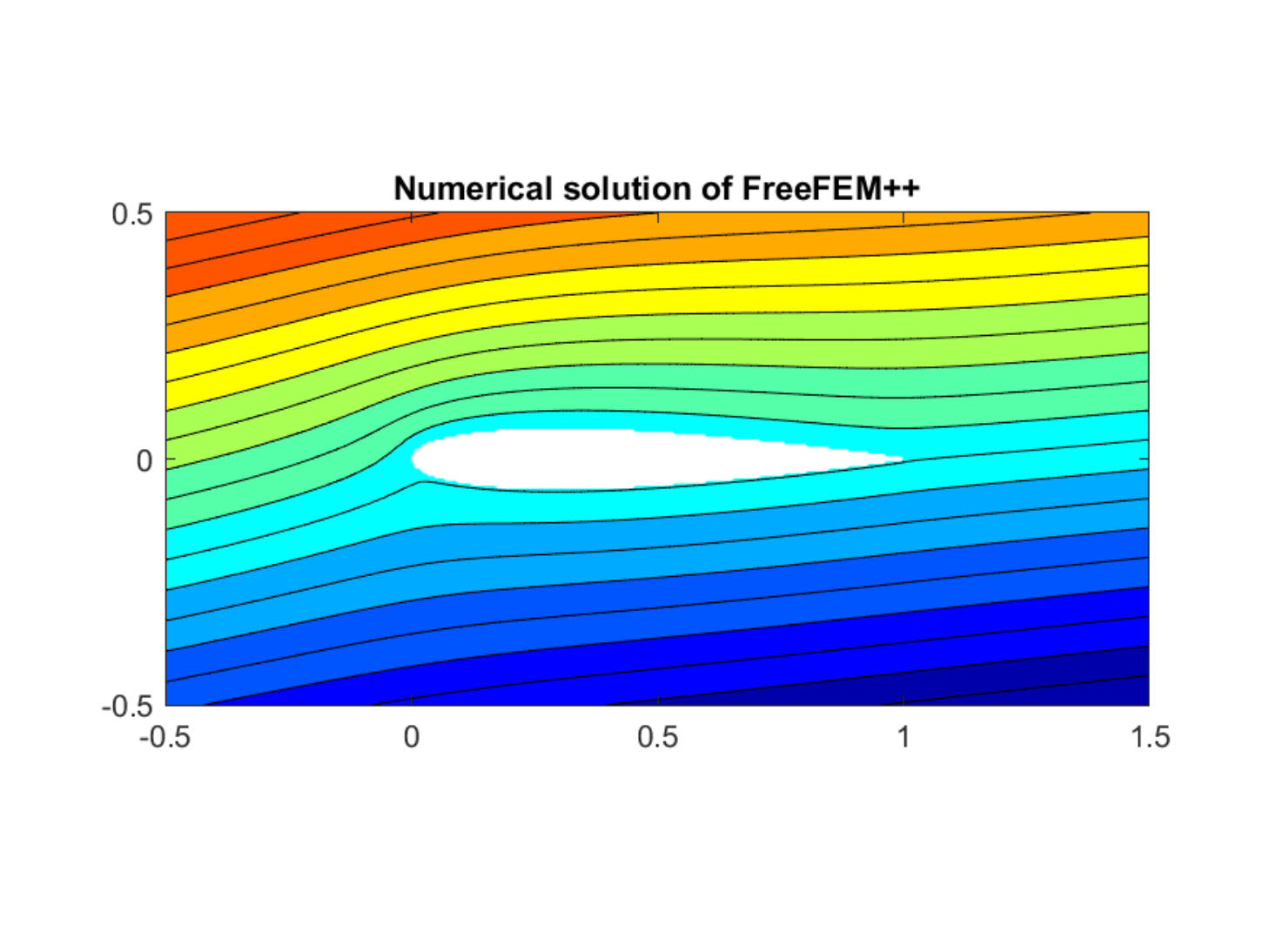}}
  \subfigure[FreeFEM solution]{\includegraphics[scale=0.4,trim=0 60 0 60,clip]{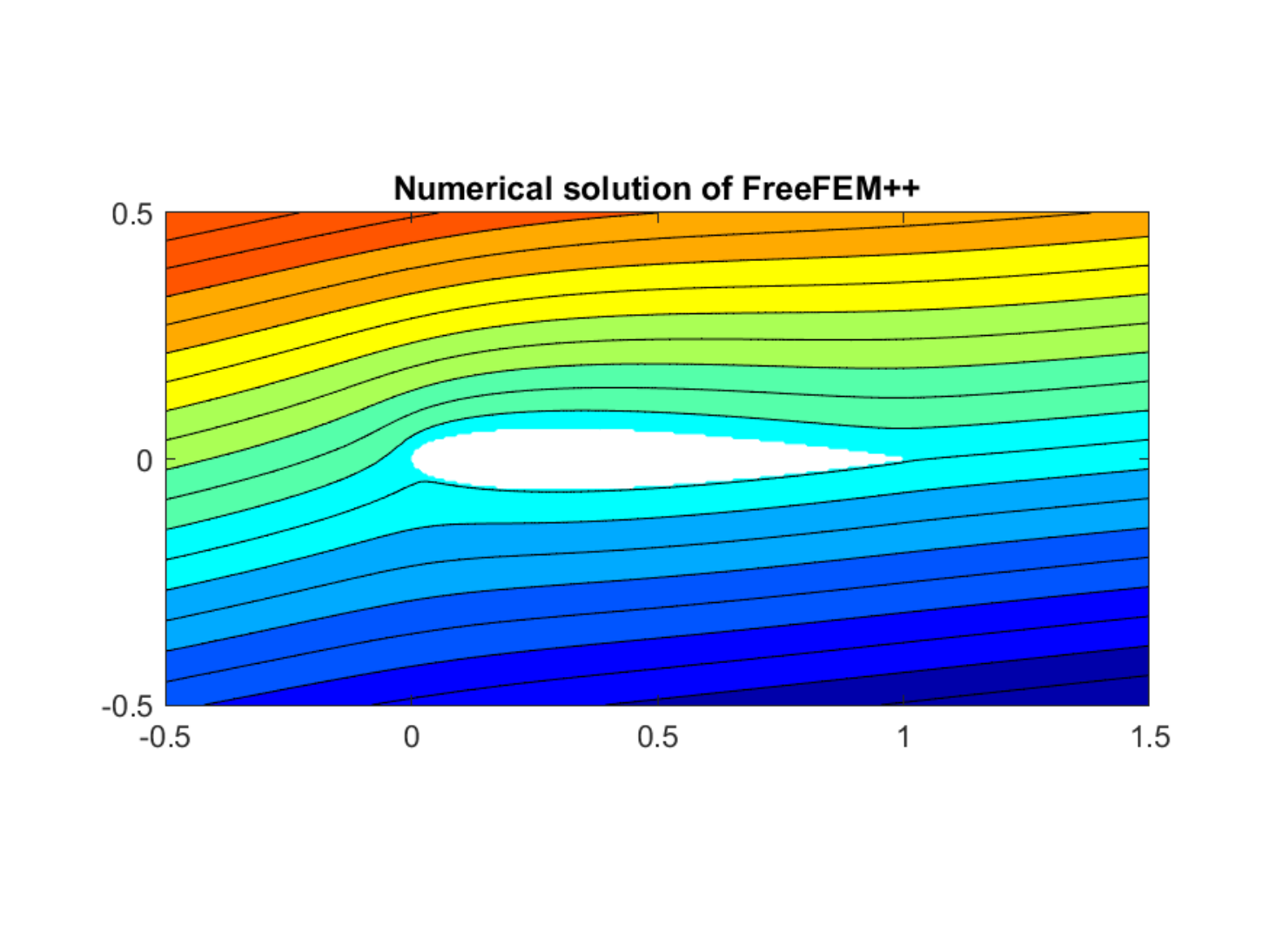}}\\
   \caption{Zoom around the NACA0012 airfoil showing the streamlines}\label{fig:naca0012_contour}
\end{figure}

We refer the reader to the FreeFEM Documentation for more details (The code is \mcode{potential.edp} in the software). The zoomed solutions of the streamlines are shown in Fig.~\ref{fig:naca0012_contour}, where the varFEM solution is well matched with the one given by FreeFEM.

\subsection{Newton method for the steady Navier-Stokes equations}

For the introduction of the problem, please refer to FreeFEM Documentation (Subsect.~2.12 - Newton Method for the Steady Navier-Stokes equations).

In each iteration, one needs to solve the following variational problem: Find $(\delta \bb{u}, \delta p)$ such that
\[DF(\delta \bb{u}, \delta p; \bb{u}, p) = F(\bb{u},p),\]
where $\bb{u}$ and $p$ are the solutions given in the last step, and
\begin{align*}
& DF(\delta \bb{u}, \delta p; \bb{u}, p) = \int_\Omega ((\delta \bb{u} \cdot \nabla)\bb{u}) \cdot \bb{v} +
(( \bb{u} \cdot \nabla) \delta\bb{u}) \cdot \bb{v} +
\nu \nabla \delta \bb{u} : \nabla v - q \nabla \cdot \delta \bb{u},  \\
& F(\bb{u}, p) = \int_\Omega ((\bb{u} \cdot \nabla)\bb{u}) \cdot \bb{v} +
\nu \nabla \bb{u} : \nabla v - p \nabla \bb{u}: \nabla \bb{v} - q \cdot \bb{u},
\end{align*}
where $(\bb{v}, q)$ are the test functions.

The finite element spaces and the quadrature rule are
\vspace{-0.8cm}
\begin{lstlisting}
Vh = {'P2','P2','P1'}; % v = [v1,v2,q] --> [v1,v2,v3]
quadOrder = 7;
vstr = {'v1','v2','q'}; ustr = {'du1','du2','dp'};
\end{lstlisting}
Here \mcode{vstr} and \mcode{ustr} are for the test and trial functions, respectively.

In the following, we only provide the detail for the assembly of the first term of $DF$, i.e.,
\begin{equation}\label{tripleNSiter}
\int_\Omega ((\delta \bb{u} \cdot \nabla)\bb{u}) \cdot \bb{v}=
\int_\Omega u_{1,x} \cdot  v_1 \delta u_1
+ u_{1,y} \cdot  v_1 \delta u_2
+ u_{2,x} \cdot v_2 \delta u_1
+ u_{2,y} \cdot v_2 \delta u_2.
\end{equation}
In the iteration, $u_{1,x}$, $u_{1,y}$, $u_{2,x}$ and $u_{2,y}$ are the known coefficients, which can be obtained from the finite element function $\bb{u}$ in the last step. Let us discuss the implementation of varFEM:
\begin{itemize}
  \item The initial data are
  \vspace{-0.8cm}
\begin{lstlisting}
% initial data
u1 = @(p) double( (p(:,1).^2 + p(:,2).^2 > 2) );
u2 = @(p) 0*p(:,1);
p = @(p) 0*p(:,1);
\end{lstlisting}
  We must convert them into the dof vectors in view of the iteration.
    \vspace{-0.8cm}
\begin{lstlisting}
uh1 = interp2d(u1,Th,Vh{1});
uh2 = interp2d(u2,Th,Vh{2});
ph = interp2d(p,Th,Vh{3});
\end{lstlisting}
  \item For numerical integration, it is preferable to provide the coefficient matrices for these coefficient functions.
      \vspace{-0.8cm}
\begin{lstlisting}
u1c = interp2dMat(uh1,'u1.val',Th,Vh{1},quadOrder);
u2c = interp2dMat(uh2,'u2.val',Th,Vh{2},quadOrder);
pc = interp2dMat(ph,'p.val',Th,Vh{3},quadOrder);

u1xc = interp2dMat(uh1,'u1.dx',Th,Vh{1},quadOrder);
u1yc = interp2dMat(uh1,'u1.dy',Th,Vh{1},quadOrder);
u2xc = interp2dMat(uh2,'u2.dx',Th,Vh{2},quadOrder);
u2yc = interp2dMat(uh2,'u2.dy',Th,Vh{2},quadOrder);
\end{lstlisting}
Here, \mcode{interp2dMat.m} provides a way to get the coefficient matrices only from the dof vectors.
  \item The triple \mcode{(Coef,Test,Trial)} for \eqref{tripleNSiter} is then given by
      \vspace{-0.8cm}
\begin{lstlisting}
Coef = { u1xc, u1yc, u2xc, u2yc};
Test  = { 'v1.val', 'v1.val', 'v2.val', 'v2.val'};
Trial = { 'du1.val', 'du2.val', 'du1.val', 'du2.val'};
\end{lstlisting}
A complete correspondence of $DF(\delta \bb{u}, \delta p; \bb{u}, p)$ can be listed in the following:
      \vspace{-0.8cm}
\begin{lstlisting}
Coef = { u1xc, u1yc, u2xc, u2yc,  ... % term 1
    u1c, u2c, u1c, u2c,  ... % term 2
    nu, nu, nu, nu, ... % term 3
    -1, -1, ... % term 4
    -1, -1, ... % term 5
    -eps ... % stablization term
    };
Test  = { 'v1.val', 'v1.val', 'v2.val', 'v2.val', ... % term 1
    'v1.val', 'v1.val', 'v2.val', 'v2.val', ... % term 2
    'v1.dx', 'v1.dy', 'v2.dx', 'v2.dy', ... % term 3
    'v1.dx', 'v2.dy', ... % term 4
    'q.val', 'q.val', ... % term 5
    'q.val' ... % stablization term
    };
Trial = { 'du1.val', 'du2.val', 'du1.val', 'du2.val',... % term 1
    'du1.dx', 'du1.dy', 'du2.dx', 'du2.dy', ... % term 2
    'du1.dx', 'du1.dy', 'du2.dx', 'du2.dy', ... % term 3
    'dp.val', 'dp.val', ... % term 4
    'du1.dx', 'du2.dy', ... % term 5
    'dp.val' ... % stablization term
    };
\end{lstlisting}

  \item The computation of the stiffness matrix reads
\vspace{-0.8cm}
\begin{lstlisting}
[Test,Trial] = getStdvarForm(vstr, Test, ustr, Trial);
[kk,info] = int2d(Th,Coef,Test,Trial,Vh,quadOrder);
\end{lstlisting}
  Here \mcode{getStdvarForm.m} transforms the user-defined notation to the standard one.
\end{itemize}

\begin{figure}[H]
  \centering
  \includegraphics[scale=0.5]{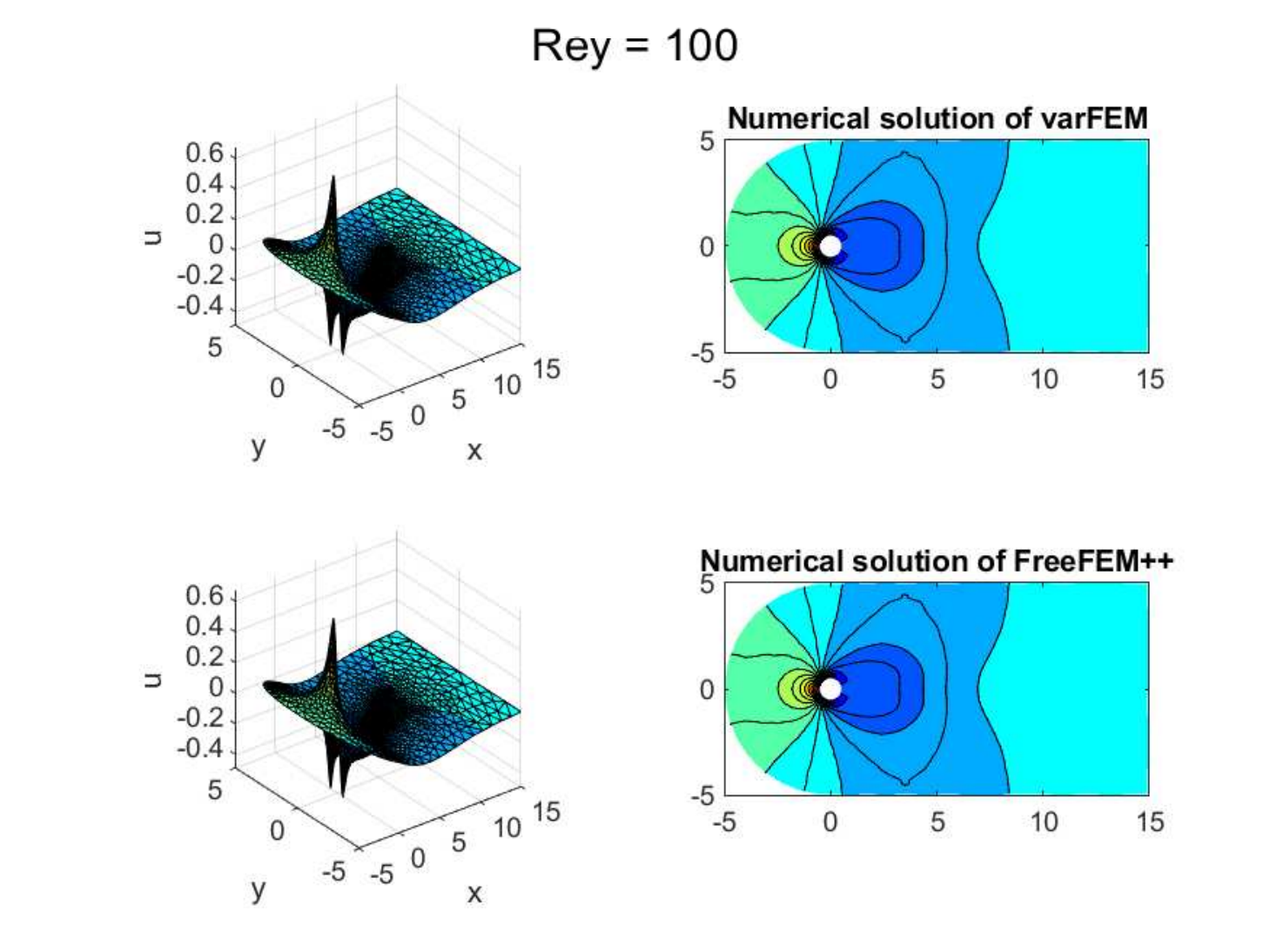}
  \includegraphics[scale=0.5]{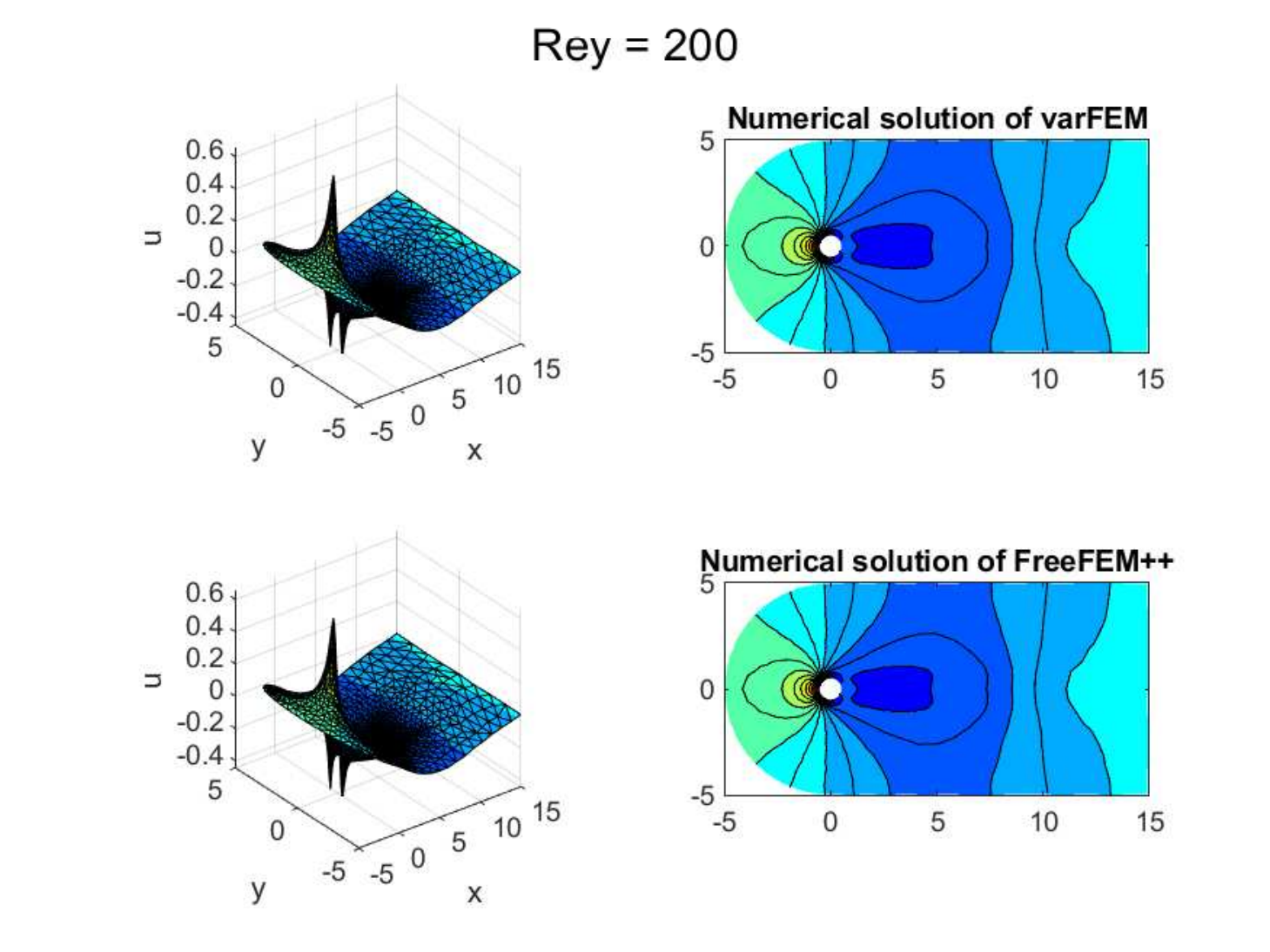}\\
   \caption{The pressure $p$ of a flow}\label{fig:NSNewton}
\end{figure}

The complete test script is available in \mcode{test\_NSNewton.m}. The varFEM solutions and the FreeFEM solutions are shown in Fig.~\ref{fig:NSNewton}.

\subsection{Optimal control}

\subsubsection{The gradient is not provided}

This is an exmple given in FreeFem Documentation: Release 4.6 (Subsection 2.15 - Optimal Control).

For a given target $u_d$, the problem is to find $u$ such that
\[\min_{z\in \mathbb{R}^3} J(z) = \int_E ( u - u_d)^2 = \int_\Omega I_E ( u - u_d)^2, \quad z = (b,c,d),\]
where $E \subset \Omega$, $I_E$ is an indication function, and $u$ is the solution of the following PDE:
\[\begin{cases}
-\nabla ( \kappa(b,c,d) \cdot \nabla u) = 0  & \mbox{in} ~~ \Omega, \\
u = u_\Gamma & \mbox{on}~~\Gamma \subset \partial \Omega.
\end{cases}\]
Let $B,C$ and $D$ be the separated subsets of $\Omega$. The coefficient $\kappa$ is defined as
\[\kappa(x) = 1 + b I_B(x) + c I_C(x) + d I_D(x), \quad x \in \Omega.\]

For fixed $z = (b,c,d)$, one can solve the PDE to obtain an approximate solution $u_h$ and the approximate objective function:
\[J_h(z,u_h,u_d) : = \int_\Omega I_E ( u_h(z) - u_d)^2.\]
We use the build-in function \mcode{fminunc.m} in Matlab to find the (local) minimizer. To this end, we first establish a function to get the PDE solution:
\vspace{-0.8cm}
\begin{lstlisting}
%% Problem for the PDE constraint
function uh = PDEcon(z,Th)

    % Parameters
    Vh = Th.Vh; quadOrder = Th.quadOrder;
    Ib = @(p) (p(:,1).^2 + p(:,2).^2 < 1.0001);
    Ic = @(p) ( (p(:,1)+3).^2 + p(:,2).^2 < 1.0001);
    Id = @(p) (p(:,1).^2 + (p(:,2)+3).^2 < 1.0001);


    % Bilinear form
    Coef = @(p) 1 + z(1)*Ib(p) + z(2)*Ic(p) + z(3)*Id(p);
    Test = 'v.grad';
    Trial = 'u.grad';
    kk = assem2d(Th,Coef,Test,Trial,Vh,quadOrder);

    % Linear form
    ff = zeros(size(kk,1),1);

    % Dirichlet boundary condition
    gD = @(p) p(:,1).^3 - p(:,2).^3;
    on = 1;
    uh = apply2d(on,Th,kk,ff,Vh,gD);
end
\end{lstlisting}

The cost function is then given by
\vspace{-0.8cm}
\begin{lstlisting}
%% Cost function
function err = J(z,ud,Th)
    Ie = @(p) ((p(:,1)-1).^2 + p(:,2).^2 <=4);
    uh = PDEcon(z,Th);

    fh = Ie(Th.node).*(uh-ud).^2;
    err = integral2d(Th,fh,Th.Vh,Th.quadOrder);
end
\end{lstlisting}

Given a vector $z_d = [2,3,4]$, we can construct an ``exact solution'' solution $u_d$ by solving the PDE. The minimizer is then given by
\vspace{-0.8cm}
\begin{lstlisting}
%% The constructed solution
zd = [2, 3, 4];
ud = PDEcon(zd,Th);

%% Find the mimimizer
options.LargeScale = 'off';
options.HessUpdate = 'bfgs';
options.Display = 'iter';
z0 = [1,1,1];
zmin = fminunc(@(z) J(z,ud,Th), z0, options)
\end{lstlisting}

In the Matlab command window, one can get the following information:
{\tiny\begin{verbatim}
                                                             First-order
      Iteration  Func-count       f(x)        Step-size       optimality
          0           4          30.9874                          77.2
          1           8          9.87548      0.0129606           21.4
          2          12          6.27654              1           12.2
          3          16          4.76889              1            5.3
          4          20          4.47092              1           3.03
          5          24          3.46888              1           5.66
          6          28          2.13699              1           5.77
          7          32          1.06649              1           4.42
          8          36         0.542858              1           3.14
          9          40         0.295327              1           2.18
         10          44         0.179846              1            1.7
         11          48         0.114221              1           1.28
         12          52         0.063295              1          0.842
         13          56        0.0242571              1          0.366
         14          60       0.00729547              1          0.258
         15          64       0.00170544              1          0.107
         16          68      3.95891e-05              1         0.0257
         17          72      3.54227e-07              1        0.00546
         18          76      3.36066e-10              1       0.000159
         19          80      5.54979e-13              1       3.88e-06

     Local minimum found.

     Optimization completed because the size of the gradient is less than
     the value of the optimality tolerance.

     <stopping criteria details>

     zmin =

         2.0000    3.0000    4.0000
\end{verbatim}}

\begin{figure}[H]
  \centering
  \includegraphics[scale=0.45]{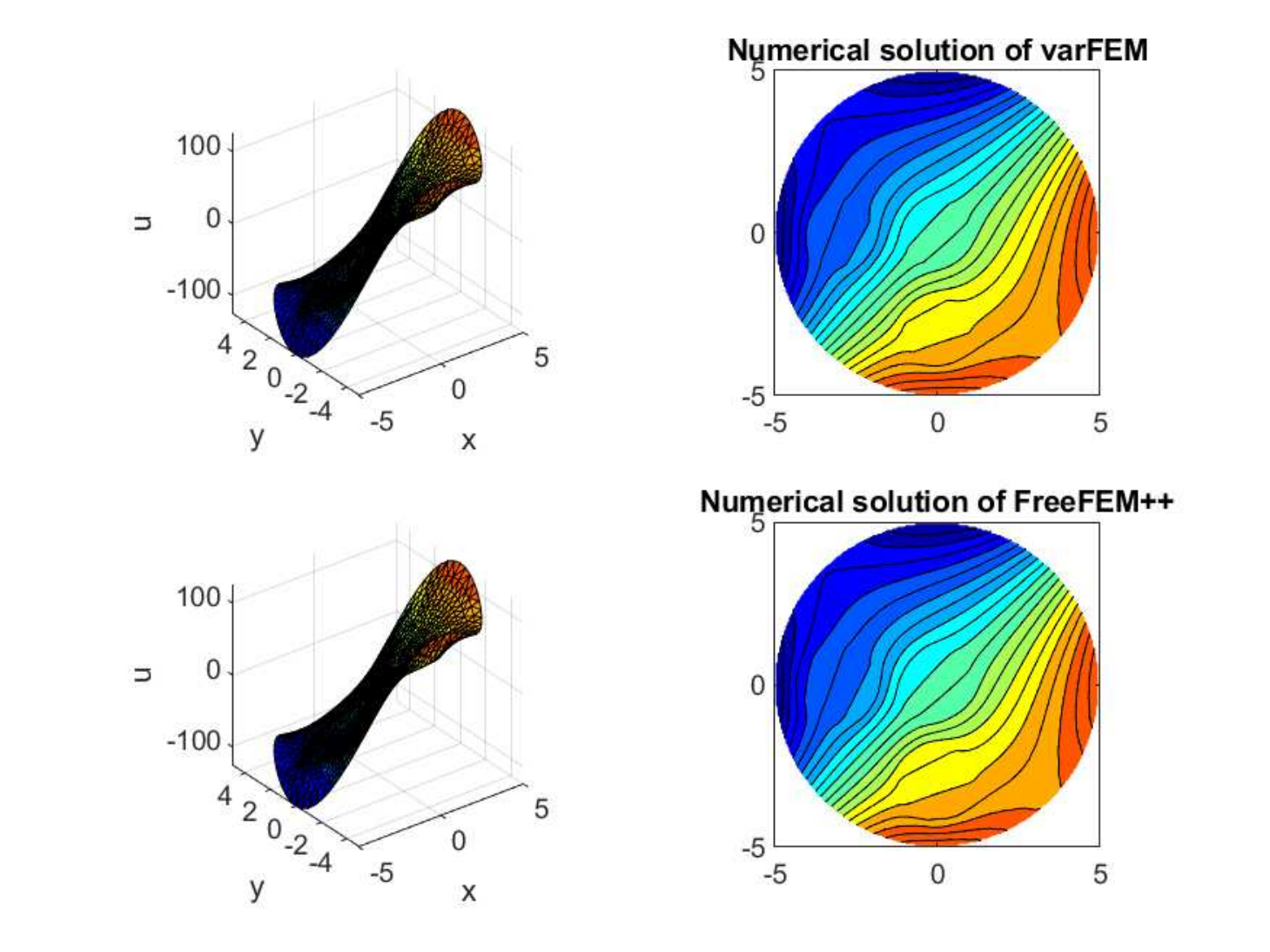}\\
  \caption{Numerical solutions of the optimal control problem}\label{fig:Optimal}
\end{figure}

We can observe that the minimizer is found after 19 iterations, with the solutions displayed in Fig.~\ref{fig:Optimal}.

\subsubsection{The gradient is provided}

For the example given in FreeFEM, the optimization problem is solved by the quasi-Newton BFGS method:
\vspace{-0.8cm}
\begin{lstlisting}
BFGS(J, DJ, z, eps=1.e-6, nbiter=15, nbiterline=20);
\end{lstlisting}
Here, \mcode{DJ} is the derivatives of $J$ with respect to $b,c,d$. We also provide an implementation when the gradient is available. In this case, the cost function should be modified as
\vspace{-0.8cm}
\begin{lstlisting}
function [err,derr] = J(z,ud,Th)
end
\end{lstlisting}
Here \mcode{err} and \mcode{derr} correspond \mcode{J} and \mcode{DJ}, respectively. For the details of the implementation, please refer to the test script \mcode{test\_optimalControlgrad.m} in varFEM. The minimizer is then captured as follows.
\vspace{-0.8cm}
\begin{lstlisting}
%% Find the mimimizer
options = optimoptions('fminunc','Display','iter',...
    'SpecifyObjectiveGradient',true); % gradient is provided
z0 = [1,1,1];
zmin = fminunc(@(z) J(z,ud,Th), z0, options)
\end{lstlisting}
Note that Line 3 indicates that the gradient is provided.

The minimizer is also found, with the printed information given as
{\tiny\begin{verbatim}
      Iteration  Func-count       f(x)        Step-size       optimality
          0           1          30.9874                          77.2
          1           2          9.87549      0.0129606           21.4
          2           3          6.27653              1           12.2
          3           4          4.76889              1            5.3
          4           5          4.47091              1           3.03
          5           6          3.46888              1           5.66
          6           7            2.137              1           5.77
          7           8          1.06648              1           4.42
          8           9         0.542849              1           3.14
          9          10         0.295333              1           2.18
         10          11         0.179846              1            1.7
         11          12         0.114223              1           1.28
         12          13        0.0632962              1          0.842
         13          14        0.0242572              1          0.366
         14          15       0.00729561              1          0.258
         15          16        0.0017055              1          0.107
         16          17      3.95888e-05              1         0.0257
         17          18      3.54191e-07              1        0.00546
         18          19      3.47621e-10              1       0.000159
         19          20      4.31601e-13              1       3.87e-06

     Local minimum found.

     Optimization completed because the size of the gradient is less than
     the value of the optimality tolerance.

     <stopping criteria details>

     zmin =

         2.0000    3.0000    4.0000
\end{verbatim}}

Compared with the previous implementation, one can find that the latter approach has fewer function calls (although we need to compute the gradient by solving a PDE problem).

\section{Concluding remarks} \label{sec:conclude}

In this paper, a Matlab software package for the finite element method was presented for various typical problems, which realizes the programming style in FreeFEM. The usage of the library, named varFEM, was demonstrated through several examples.
Possible extensions of this library that are of interest include three-dimensional problems and other types of finite elements.

%%\newpage
%%\section*{References}
%
%\bibliographystyle{plain} %plain, unsrt
%\bibliography{Refs}

\begin{thebibliography}{1}

\bibitem{ChenL-iFEM-2009}
L.~Chen.
\newblock {iFEM}: an integrated finite element method package in {MATLAB}.
\newblock Technical report, University of California at Irvine, 2009.

\bibitem{FreeFEM}
F.~Hecht.
\newblock New development in freefem++.
\newblock {\em J. Numer. Math.}, 20(3-4):251--265, 2012.

\end{thebibliography}

%\end{CJK*}

\end{document}